\documentclass{article}

\usepackage{authblk}
\usepackage[dvipsnames]{xcolor}
\usepackage[12pt]{extsizes}
\usepackage[utf8]{inputenc}
\usepackage[english]{babel}
\usepackage{fancybox}		
\usepackage{amsmath}
\usepackage{amsthm}
\usepackage{comment}
\usepackage{amsfonts}
\usepackage{enumerate}
\usepackage[shortlabels]{enumitem}
\usepackage{graphicx}
\usepackage[left=2cm, right =2cm, top = 2cm, bottom =2cm]{geometry}
\usepackage{fullpage}
\usepackage{hyperref}
\usepackage{xcolor}
\usepackage{minted}
\usepackage{stmaryrd} 

\usepackage{tikz}
\usepackage{tikz-qtree}
\usepackage{tikz-cd}
\usepackage{graphicx}
\usepackage{array}
\usepackage{mathabx}
\usepackage{amssymb}

\usepackage{mathtools}
\makeindex

\newcommand{\Lc}{\mathbb{L}}

\newcommand{\Q}{\mathbb{Q}}

\newcommand{\Qb}{\overline{\Q}}
\newcommand{\Qpb}{\overline{\Q_p}}

\newcommand{\Z}{\mathbb{Z}}

\newcommand{\F}{\mathbb{F}}

\newcommand{\N}{\mathbb{N}}

\newcommand{\Uset}{\mathcal{U}}

\newcommand{\HJ}{\Hop^1 ( \Q ,J[2])}
\newcommand{\HJk}{\Hop^1(k,J[2])}
\newcommand{\Sel}{\Selmer^{(2)}}
\newcommand{\muT}{\mu_2^{\OTC}}
\newcommand{\muDelta}{\mu_2^\Delta}
\newcommand{\mT}{\muT / \mu_2}
\newcommand{\mDelta}{\muDelta/\mu_2}
\newcommand{\FT}{\F_2^{\OTC}}
\newcommand{\FDelta}{\F_2^\Delta}
\newcommand{\Hmu}{\Hop^1(\Q,\mT)}

\newcommand{\HDelta}{\Hop^1(\Q,\mDelta)}
\newcommand{\HDeltak}{\Hop^1(k,\mDelta)}

\newcommand{\etale}{L^*/L^{*,2}\Q^*}
\newcommand{\etalek}{L_k^*/L_k^{*,2}k^*}

\newcommand{\qs}{q_{symp}}

\newcommand{\OTC}{\Delta_{Odd}}
\newcommand{\TC}{\Delta_{TC}}
\newcommand{\ETC}{\Delta_{Even}}
\newcommand{\kb}{\overline{k}}

\newcommand{\JQ}{J(\Q)/2J(\Q)}
\newcommand{\JQp}{J(\Q_p)/2J(\Q_p)}
\newcommand{\Jk}{J(k)/2J(k)}
\newcommand{\E}{E^\vee}
\newcommand{\Rw}{R^\vee}
\newcommand{\kF}{\mathcal{K}}
\newcommand{\sKF}{\widetilde{\kF}}
\newcommand{\Fg}{\F_2^{2g}}

\newcommand{\Sp}{\SpOp_{2g}(\F_2)}

\newcommand{\SelF}{\Selmer_{Fake}}
\newcommand{\SelT}{S_\mathcal{T}}
\newcommand{\wSel}{\Selmer_{Fake,Inter}}
\newcommand{\HE}{\Hop^1(k,\E)}
\newcommand{\HR}{\Hop^1(k,\Rw)}

\newcommand{\piOTC}{\pi_{\OTC}}

\newcommand{\Oc}{\mathcal{O}}
\newcommand{\LS}{L(2,\Sset)}

\newcommand{\sLS}{\widetilde{\LS}}

\newcommand{\Sset}{\mathcal{S}}
\newcommand{\Tset}{\mathcal{T}}
\newcommand{\W}{W^\wedge}

\newcommand{\rpztors}{\rho_{2}}
\newcommand{\rpzOTC}{\rho_{odd}}
\newcommand{\rpzTC}{\rho_{TC}}
\newcommand{\rpzF}{\rho_\F}

\newcommand{\PK}{(\Pic_X/\langle K_0 \rangle) [2]}
\newcommand{\<}{\langle}
\newcommand{\ran}{\rangle}

\newcommand{\ra}{\rightarrow}
\newcommand{\hra}{\hookrightarrow}
\newcommand{\surj}{\twoheadrightarrow}

\usepackage[OT2,T1]{fontenc}
\DeclareSymbolFont{cyrletters}{OT2}{wncyr}{m}{n}
\DeclareMathSymbol{\Sha}{\mathalpha}{cyrletters}{"58}

\newcommand{\divi}{\operatorname{div}}

\def\res{\mathrm{res}\,}

\newtheorem{theo}{Theorem}[section]
\newtheorem{de}[theo]{Definition}
\newtheorem{prop}[theo]{Proposition}
\newtheorem{cor}[theo]{Corollary}

\theoremstyle{remark}
\newtheorem{rem}[theo]{Remark}
\theoremstyle{remark}
\newtheorem{example}[theo]{Example}

\newcommand{\Pic}{\operatorname{Pic}}
\newcommand{\Jac}{\operatorname{Jac}}
\newcommand{\Selmer}{\operatorname{Sel}}
\newcommand{\Hop}{\operatorname{H}}

\newcommand{\DivOp}{\operatorname{Div}}
\newcommand{\rank}{\operatorname{rank}}
\newcommand{\Arf}{\operatorname{Arf}}
\newcommand{\Aut}{\operatorname{Aut}}
\newcommand{\coker}{\operatorname{coker}}
\newcommand{\im}{\operatorname{Im}}

\newcommand{\Frob}{\operatorname{Frob}}

\newcommand{\SpOp}{\operatorname{Sp}}

\newcommand{\bs}{\bigskip}

\renewcommand{\contentsname}{Sommaire} 
\everymath{\displaystyle}
\usepackage{url}

\usepackage{algpseudocode}
\usepackage{algorithm}

\newcommand{\Nc}{\mathcal{N}}
\newcommand{\NQ}[1]{(\Q^*/\Q^{*,2})^{#1}/\Q^*}
\newcommand{\NQo}{\NQ{w_o}}
\newcommand{\NQe}{(\Q^*/\Q^{*,2})^{w_e}}
\newcommand{\ev}{\operatorname{ev}}

\usepackage{framed}

\begin{document}

\title{\vspace{-10mm}Explicit generalized $2$-descent on Jacobians of non-hyperelliptic curves of moderate genus}
\author{Thibaut Misme\thanks{\href{mailto:mismet@tcd.ie}{mismet@tcd.ie}}}
\affil{\scriptsize{Trinity College Dublin (TCD)}}
\maketitle
\vspace{-5mm}

\begin{abstract}
	We implement and present an automatized algorithm performing 2-descent on Jacobians of general curves. It follows the method of Bruin-Poonen-Stoll based on Odd Theta Characteristics. We also extend our algorithm to perform 2-descent using all Theta Characteristics. We call it TC-descent. 
    After reviewing generalized descent on non-hyperelliptic curves, we run explicit TC-descents on Jacobians of curves from their equations only, ignoring specific features these curves may have. In particular, we achieve a successful 2-descent on a non-hyperelliptic genus 5 curve with splitting TC-orbits.
\end{abstract}

\renewcommand{\contentsname}{Contents}
\tableofcontents 






\newpage

    
\section*{Introduction}
\addcontentsline{toc}{section}{Introduction}

    \subsection*{Context}
	Computing the set of rational points of a given curve is a very old problem. Nowadays one standard way stands out (See \cite{StollRationalPoints}, \cite{BruinRationalPoints} and \cite{PoonenBook}): First, ask a computer to find as many rational points as it can using brute-force based techniques; second, prove that these points are the only ones via Chabauty-like methods (See \cite{Coleman}); see \cite{ChabautyColemanSurvey} for a survey. But in order to be able to apply this process, we need to ensure that the Mordell-Weil rank of the Jacobian is (strictly) lower than the genus of the curve. Very efficient and explicit algorithms, namely $p$-descent algorithms, have been developed for elliptic curves focusing on Galois Cohomology arising from $p$-torsion (See \cite{CEbible}). The case $p=2$ is the easiest, most common, and most computationally efficient one; furthermore, it naturally extends to hyperelliptic curves (See \cite{StollCH}).\par
    
    The $2$-descent method has been very much studied and is now quite well understood on the theoretical side. More precisely, it aims to bound the rank of the Jacobian by computing cardinals of Selmer groups. When it comes to Jacobians of general curves, Bruin, Poonen and Stoll (BPS) offered a very general theoretical framework to perform generalized descent method in \cite{BPS} based on Odd Theta Characteristics (OTC). They present explicit examples for plane quartics; Fu and Stoll perform the same process on a non-hyperelliptic curve of genus $4$ in \cite{StollFu}; Bruin and Lewis collect statistics by running many examples over plane quartics whose bitangents are all rational in \cite{BruinRatBitangents}. 
    
    \subsection*{Present work}
    In opposite to the hyperelliptic case, automatized 2-descent has remained intractable in practice on Jacobians of general curves of genus strictly greater than $4$. It is substantially due to the lack of means to compute the Galois action of the set of OTC.  
    Building upon Mascot’s algorithm \cite{Mascot}, the author has designed an algorithm \cite{TitiOTC} precisely to fulfil this gap. In the present article, we exploit these techniques to implement automatized 2-descent in the meaning of BPS (See \cite{BPS}) on non-hyperelliptic curves.\par
    Contrary to the methods used in \cite{BPS}, our algorithm \cite{TitiOTC} does not rely on the geometrical counterpart of the OTC, such as bitangents and tritangents for genera 3 and 4, to compute the Galois action on them. It allows us to reach more invariants but more importantly to access the whole set of Theta Characteristics (TC). This is the main novelty of our approach.\par

    We take advantage of this computational asset to perform a so-called Fake descent with all TC. We call it \textit{TC-descent}. It always gives descent maps whose kernels are smaller than with OTC-descent as presented in \cite{BPS}.
    We study, both theoretically and computationally, the kernels of the TC-descent map in a work to appear. In particular, our approach to obtain the TC makes the computation of the kernels easier. We also state theoretical and practical conditions for injectivity. See section \ref{section explicit Galois} of the present article for more details.
    Note that TC-descent may also be interesting for odd genus hyperelliptic curve whose descent kernel control is left over by the classic article \cite{StollCH} from Stoll.\par
    
    We present several explicit TC-descents ran by the algorithms implemented by the author.
    They apply to general curves, and do not rely on the geometry nor special features of the curves under consideration. Nevertheless, where results had previously been obtained by more specific methods for particular curves, invoking specific properties that these curves enjoy, we are able to rederive these results without relying on these specific properties anymore.
    We also manage a successful TC-descent on a genus 5 curve with splitting TC-orbits, which had never been done before.

	\subsection*{Inputs of the algorithm}
	We present an algorithm performing Fake descent in the meaning of BPS (See \cite{BPS}) using (Odd) Theta Characteristics. It asks for
    \begin{enumerate}
    \item A single equation $f(x,y)\in \Q[x,y]$ whose associated (singular) projective model in $\mathbb{P}^2$ is birational to the curve $X$ we are interested in.
    \item A rational point $P_0$ on the model $f$ (or a degree 1 rational divisor).
    \item A description of the Galois group $G$ of $\OTC$ or $J[2]$ as embedded in $\Sp$.
\end{enumerate}
    
	We can try to search input 2. by brute-force if not provided. Methods to obtain input 3. from the data we compute anyway are described in \cite{BPS}. When successful, the algorithm returns the dimension of the rational 2-torsion and a bound for the size of the $2$-Selmer group. It therefore also produces an upper bound for the rank of the Jacobian.\par
    The main drawback of our method is the requirement of a rational point, as many curves do not have any, and checking its existence is a famously hard question in Diophantine geometry.
	
	\subsection*{Organisation of the article}
	In the first chapter, we review what 2-descent is in a way that eases the understanding of Fake descent \`a la BPS.\\
	In chapter 2, we mainly remind Fake descent as it is presented in \cite{BPS}.\\
	In chapter 3, we expose our implementation of descent in a computational point of view.\\
	In chapter 4, we present TC-descent.\\ 
	In chapter 5, we expose a few ways to certify the rank once an upper bound has been found.\\
	In chapter 6 and 7, we display the examples that the author ran. To assess trust and correctness of the method, we recover several examples of the literature using only the generic inputs that the algorithm of the author demands.\\
    Chapter 7 is dedicated to a TC-descent on a genus 5 curve with splitting TC orbits.\par
	
	Nothing is new in the two first chapters unless explicitly mentioned. Chapter 3 and 5 relies in part on the common literature of 2-descent. Chapters 4, 6 and 7 are entirely the work of the author.
	
	\subsection*{Acknowledgements}
	I thank my PhD supervisor Nicolas Mascot for continuous valuable discussions and advice about the present article. This work has been accomplished at Trinity College Dublin (TCD), but was written during the author's visit at the Aix-Marseille University (amU). The algorithm has been implemented in the free and open-source package PARI/GP \cite{PARI2}.\bs
    
	The author would be happy to try performing 2-descent on other non-hyperelliptic curves than presented in the examples section. The reader should not hesitate to contact him with a given curve (of moderate genus) to perform descent on along with a short motivation of the original purpose.\bs

    \noindent
	\textbf{Keywords}: 2-descent, Fake descent, Theta Characteristic, Mordell-Weil rank.\\
    \textbf{Notations}: 
	In the present article, by \textit{nice curve over a field $k$} we mean a curve that is smooth, projective and geometrically integral. 

\section{Toward generalized 2-descent on Jacobians of curves}

We review in this section some known facts related to 2-descent in a form that we want as pedagogical as possible rather than formal; \cite{BPS} gives a full description. First, we give an introduction to what 2-descent is in an interpretation that eases the understanding of the 2-descent approach of Bruin, Poonen and Stoll (BPS) \cite{BPS}. After exposing some short useful facts about (Odd) Theta Characteristics we will present the BPS approach in the next chapter.

\subsection{Classic 2-descent}

In its modern, computational point of view, 2-descent is a method that aims to bound the Mordell-Weil rank of the Jacobian of a smooth curve. It is often a substantial step in the task of determining the set of rational points on the curve; see \cite{StollRationalPoints}, \cite{BruinRationalPoints} for surveys on computing rational points; see \cite{PoonenBook} for a thorough coverage of the topic. See \cite{CEbible} for a complete approach on $n$-descent over elliptic curves; \cite{StollECarticle} for an efficient survey. One can present Selmer groups under a more geometrical point of view, namely as parametrizing twists of the curves with specific features; see \cite{StollEC}, \cite{GeometricDescent} and \cite{GeometricDescentCH}. In the present article we will only tackle 2-descent in a computational setting. For hyperelliptic curves, see \cite{StollCH}; for $p$-covers of $\mathbb{P}_1$, see \cite{PoonenSchaefer}.\bs

The 2-descent method is based on the Mordell-Weil theorem:
\begin{theo}[Mordell-Weil] \mbox{}\\
    Let $K$ be a number field and $A$ an abelian variety over $K$.\\
    The group $A(K)$ of $K$-rational points of $A$ is finitely generated. In other words, there exists an integer $r\in \N$ (called the Mordell-Weil rank of $A/K$) usually denoted as $\rank_K(A)$ and a \text{\normalfont{finite}} group $A_{tors}$ such that we have the isomorphism of groups
    $$ A(K) \simeq A_{tors} \times \Z^r.$$
\end{theo}

After a heavy use of height theory (which is actually the main part), the proof breaks down to show that the group $A(K)/mA(K)$ is finite for any $m\in \N$ (See \cite{DiophantineGeometry} Part C). This is achieved by embedding it into the infinite Galois cohomology group $\Hop^1(K,A[m])$ derived from the $m$-torsion points $A[m]$ and cutting out by local conditions to obtain a finite group called the \textit{Selmer group}. The Selmer group contains $A(K)/mA(k)$ by design, and its finiteness completes the proof. In fact, the cardinal of the Selmer group provides an upper bound on the \textit{Mordell-Weil rank} $r$ of $A$ as consequently $A(K)/mA(K) \simeq A[m](K) \times (\Z/m\Z)^r$.\bs

The 2-descent method takes $m=2$ for computational efficiency and aims at making this part of the proof explicit. The goal is henceforth to compute the 2-Selmer group using the very techniques that arose from the proof of the Mordell-Weil theorem. More explicitly, the 2-Selmer group is defined as
$$ \Sel_K(A):= \{ \zeta \in \Hop^1(K,A[2])~|~\forall p \text{ prime, }\res_p(\zeta)\in \gamma_p(A(K_p)/2A(K_p)) \}$$ 
where $\res_p : \Hop^1(K,A[2]) \ra \Hop^1(K_p,A[2])$ is the restriction map to a decomposition group at the prime $p$ and $\gamma_p : A(K_p)/2A(K_p) \ra \Hop^1(K,A[2])$ comes from the boundary cohomology map constructed from the short exact sequence of Galois-modules $ 0 \ra A[2] \ra A \overset{2}{\ra} A \ra 0$. \par
Equivalently, the Selmer group can be defined as an intersection of inverse images: 
$$\Sel_K(A) := \bigcap_p ~ \res_p^{-1}(\gamma_p(A(K_p)/2A(K_p))). $$ The following commutative diagram unravels its definition.
\begin{equation*}
        \begin{tikzcd}
            & \Sel_K(A) \arrow[d, phantom, sloped, "\subset"]\\
            A(K)/2A(K) \arrow[r, hook, "\gamma_K"] \arrow[d,"\res_p"] & \Hop^1(K,A[2]) \arrow[d,"\res_p"]\\
            A(K_p)/2A(K_p) \arrow[r, hook,"\gamma_p"] & \Hop^1(K_p,A[2])
        \end{tikzcd}
\end{equation*}

We thus obtain a bound $\rank_K(A) \leq \dim_{\F_2} \Sel_K(A) - \dim_{\F_2} A[2](K)$, which is often sharp. Actually the error is controlled by a space which is in general out of computational reach, namely the 2-torsion of the Tate-Shafarevich group $\Sha_\Q^2(A)$. In practice the latter group is often (but not always) trivial, making the whole method quite relevant.

\subsection{An introduction to generalized 2-descent}

The need of an ambient space, where the intersection over all primes $p$ takes place in the hope of capturing $A(K)/2A(K)$ justifies the use of Galois cohomology though $\Hop^1(K,A[2])$. But as the space $\Hop^1(K,A[2])$ is not explicit, the 2-descent method projects it (via the so-called \textit{cohomology map} which we shall introduce soon) onto a more computationally-friendly space $\Hop^1(K,V)$ for a certain Galois-module $V$. The method generalized in \cite{BPS} then consists in approaching the equivalent \textit{Fake Selmer group} in $\Hop^1(K,V)$ cut out by the same local conditions than the classic Selmer group. The natural definition of the Fake Selmer group by local conditions (in finite vector spaces that you can check by brute-force) makes its computation amenable. Ideally, we would like the cohomology map $\Hop^1(\Q,A[2]) \longrightarrow \Hop^1(\Q,V)$  to be injective on the Selmer; but unfortunately, this is not always the case, and the success of the method then relies on our capacity to control the intersection of its kernel with the classic Selmer group. This can be quite a sensitive step in practice as it requires in general the computation of the Galois group of the 2-torsion $J[2]$ of the Jacobian. The next section focuses on giving more details to this part of the 2-descent process.\par
For instance, in the simplest case where $A$ is an elliptic curve $E : y^2=f(x)$ defined over $K=\Q$, classic $2$-descent can be formulated in this framework as taking $V=\mu_2^\Delta$, where $\Delta$ is the set formed by the 3 Weierstrass points attached to the roots of $f$. This gives a map $\Hop^1(\Q,E[2]) \ra \Hop^1(\Q,\mu_2^\Delta) \simeq L^*/L^{*,2}$ where $L$ is the \'etale algebra of $\Delta$. In that convenient case, the cohomology map is injective and the classic Selmer group can be recovered as the kernel of an explicit map inside the Fake Selmer group making the process much easier and more accurate than in general.\par

The Mordell-Weil theorem, along with the method's strategy, also apply to Jacobians of curves defined over $K=\Q$. Figuring out the rank of the Jacobian turns out to be an essential step for determining the rational points of the underlying curve as well when the curve is not elliptic; for instance in order to use techniques \a la Chabauty-Coleman \cite{Chabauty}.\par
The 2-descent method for elliptic curves has been the topic of an extensive study both theoretically than computationally. It generalizes quite well in its spirit (but not as easily in its practicability) for hyperelliptic curves (See \cite{StollCH}) or to $p$-covers of $\mathbb{P}_1$ (See \cite{PoonenSchaefer}). However, extending it to non-hyperelliptic curves following \cite{BPS}, which is the purpose of the present article, is more challenging. Bruin, Poonen and Stoll consider the Galois-set $\Delta:=\OTC$ formed by the Odd Theta Characteristics of the curve to base the cohomology map $\HJ \ra \Hmu$ on. We dive more into this construction for Jacobians in the next section, but first, as they play an essential role in what follows, we quickly review Theta Characteristics.

\subsection{Theta Characteristics}\label{section TC}

For a reminder about Theta Characteristics tailored to the use of the present article, see \cite[section 2.2.1.]{TitiOTC}. For a reminder about the Galois action over Theta Characteristics, see \cite[section 2.3.]{TitiOTC}. In particular, the proofs of the following facts can be found in the mentioned article.\bs

Let $k$ be a field with $char(k) \neq 2$. Fix a canonical divisor $K_X$ on $X$ such that $K_X$ is defined over $k$.

\begin{de}[Theta Characteristic] \mbox{}\\
    A Theta Characteristic (TC) over k is a divisor $\beta \in \DivOp_X(k)$ such that
    $$ 2 \beta \sim K_X.$$
    We define $\Delta(k) \subset \Pic_X (k) $ to be the set of the \normalfont{classes} \textit{of TC over $k$.\\
    A Theta Characteristic $\beta$ is called Odd (OTC) if $h^0(\beta):=\dim \Hop^0(\beta)$ is odd. It only depends on the class $[\beta]$ and therefore we set $\OTC(k) \subset \Delta(k)$ to be the} \normalfont{classes} \textit{of OTC over $k$.}\textit{ Otherwise we say that $\beta$ is even (ETC) and we set $\ETC(k) \subset \Delta(k)$ to be the} \normalfont{classes} \textit{of ETC over $k$.}
\end{de}
We recall that these sets are finite. Over an algebraic closed field,
$$ |\Delta(\kb)| = 2^{2g} \hspace{10mm} \text{and} \hspace{10mm} |\OTC(\kb)| = 2^{(g-1)}(2^g -1).$$
They are also stable under Galois. In particular, $\OTC$, $\ETC$ and $\TC$ are finite non-empty Galois sets.\par
A trivial property but of crucial importance is that any difference of TC in the Picard group gives a $2$-torsion point on the Jacobian. Fixing any TC, one also notices that any $2$-torsion point can be obtained that way. Actually, we can even restrict to OTC.
\begin{prop}\label{OTC diff} 
For every point $P\in J[2]$, there exists $\beta_1,\beta_2 \in \OTC$ such that
    $$ \beta_1 - \beta_2 = P.$$
\end{prop}

\begin{proof}
    See \cite[Cor 5.3.a.]{BPS}.
\end{proof}

Adding to the classic Galois representation $\rpztors : G_k \ra \Sp$ afforded by the $2$-torsion $J[2]$, we can also construct permutation representations derived from the natural Galois action on $\OTC$ and $\TC$:  
$$ \rpzOTC : G_k \ra \operatorname{Sym}(\OTC) \hspace{10mm} \rpzTC : G_k \ra \operatorname{Sym}(\TC).$$
Both images can be represented non-canonically by a subgroup of permutations as the underlaying Galois-sets are finite. It turns out that the fields of definition of $J[2]$, $\TC$ and $\OTC$ are all actually the same field (See \cite[Prop 2.17.]{TitiOTC}). We therefore denote it by $\Lc$. Consequently, $\OTC$, $\TC$ and $J[2]$ have the same Galois group and, equivalently, the 3 corresponding representations $\rpzOTC$, $\rpzTC$, and $\rpztors$ have images which are isomorphic as groups:
$$ \im(\rpzOTC) \simeq \im(\rpzTC) \simeq \im(\rpztors) \hra \Sp.$$
More generally, the Galois action on the OTC is determined by the one on the $2$-torsion or the TC and vice-versa.
More details can be found in \cite{SHA}.\bs

Finally, the equality between the definition fields also provides knowledge on the ramification primes of $\Lc_{Odd}=:\Lc$ by the Neron-Ogg-Shafarevich theorem: 
\begin{prop}\label{LOTC ramification}
    The extension $\Lc/k$ is unramified away from the union of the set of the primes of bad reduction of the Jacobian and the set of the primes above 2.
\end{prop}


We have an equivalent notion of unramification for cocycles. Let $M$ be a Galois-module over a nonarchimedean local field $k_p$. We denote $I_p\subset G_k$ its inertia group.
\begin{de}\label{def unramified}
    We set $\Hop^1_{unr}(k,M):= \ker(\Hop^1(k_p,M) \ra \Hop^1(I_v,M))$. If $\xi\in \Hop^1_{unr}(k,M)$, we say that the cocycle $\xi$ is unramified.
\end{de}

\section{Generalized explicit descent in the Bruin-Poonen-Stoll meaning}

We give a very short sum-up of \cite{BPS}. Also See \cite[sections 7. and 6.]{BruinRationalPoints}. Nothing in this chapter is new, and the reader should refer to \cite{BPS} for the proofs or the global picture and at the following sections for details about how the author performs explicit computations. We will try to follow the notations from the original article. Consider a curve $X$ of genus $g \geqslant 1$ defined over $\Q$. The theory straightforwardly generalises to number fields but it makes the computations much involved. We write $J := \Jac_\Q(X)$ for its Jacobian, we fix a canonical divisor $K_X$ defined over $\Q$ and we recall that the set $\OTC := \lbrace $ classes of OTC $\rbrace$ is of cardinal $2^{(g-1)}(2^g -1)$. For the rest of the article, $k$ will denote a global field or its completion, mainly $\Q$ or $\Q_p$.\bs

As we said, in order to generalize the classic 2-descent method on elliptic curves to Jacobians of arbitrary curves, we need to remedy the lack of computational effectiveness of the Galois cohomology group $\HJk$. In the edifying article \cite{BPS}, N.Bruin, B.Poonen and M.Stoll (BPS) project $\HJk$ onto a more handy cohomology group $\Hop^1(k,V)$ using morphisms $\alpha_k: \HJk \ra\Hop^1(k,V)$ of prime interest in the present article.  We will call it the \textit{(global/local) cohomology map} and describe it in the first section of this chapter.\par
The morphisms commute with the restriction to the decomposition subgroups. The latter property is essential due to the local-global colour of the Fake Selmer's description the computation is based on. We will introduce it formally in the second section.\par
BPS then deal with controlling the kernel of the cohomology map, referred to as \textit{cohomology kernel} (See Definition \ref{def CohoKernel} in the present article. We will recall how they do so in the third section (See section \ref{sectionKernel}).

\subsection{The cohomology map $\alpha_k$ and the descent map $C_k$}\label{section Descent Cohomology}

We follow the approach of \cite{BPS}, generalizing the elliptic curve case, to project $\HJk$ onto an explicit space.\bs

For any Galois-set $\Delta$, we denote $\FDelta$ the permutation-module of $\Delta$ over $\F_2$. It follows that $\FT$ is naturally equipped with a Galois action which permutes the elements of the canonical basis $(e_\beta)_{\beta \in \Delta}$ according to the elements in $\Delta$ it represents. Let $E:=(\F_2^{\Delta})_0$ be the vector-subspace of $\FT$ formed by the elements with an even number of non-zero coordinates (i.e. of trace 0).\par
BPS start from a Galois-set $\Delta$ equipped with a \underline{surjective} morphism of Galois-modules $ \pi_\Delta: V \surj J[2]$, where $V:=\F_2^\Delta$ (\textit{True Descent}) or $V:=E$ (\textit{Fake Descent}). See Definition 6.2 and Definition 6.7 in \cite{BPS} for the actual definition. In general, for a smooth curve, we can consider a Fake Descent setting by taking $\Delta:=\OTC$ to be the set of classes of Odd Theta Characteristics. Indeed, we can identify $J[2]$ with its embedding in $\PK$ to define the following map:
\begin{align*}
    \piOTC : E:=(\F_2^{\OTC})_0 & \surj J[2]\\
    \sum_{\beta \in \OTC} \mu_\beta e_\beta & \mapsto \sum_{\beta \in \OTC} \mu_\beta \beta.
\end{align*}
This is a morphism of Galois-modules, which is surjective by Proposition \ref{OTC diff}. More generally, we could take $\Delta$ to be the union of any chosen orbits of TC, as long as $\Delta$ is big enough to span $J[2]$ under the equivalent morphism $\pi_\Delta$. Under the assumption of the existence of a rational point $P_0\in X(\Q)$ on the curve, these $\Delta$ define indeed an actual \textit{Fake Descent setup} as required in \cite[Def 6.7.]{BPS}. See the section \ref{section Fake Descent setup} of the present article for more details.
\begin{rem}
Notice that the morphism is well-defined over $E$ but could not be defined over the whole space $\FT$ in general, thus preventing a True Descent. In contrast, an elliptic curve $E$ admits a True Descent setting by taking $\Delta:= E[2]\setminus 0$. This is actually a general fact: every curve with a rational OTC $\beta_0$ admits a True Descent setting by taking $\Delta:= \OTC\setminus\{ \beta_0 \}$. By slightly modifying the morphism $\piOTC$, one can equip $\OTC$ with a True Descent setting even if the rational TC is even.\bs
\end{rem}

Let $R := \ker(E \rightarrow J[2])$ be the kernel of the morphism $\pi_\Delta$. It can be shown that it is induced by the Syzygetics quadruples the same way the pairs of OTC induce $E$. By dualising the short exact sequence of Galois-modules
$$ 0 \ra R \ra E \ra J[2] \ra 0$$
and noticing that $J[2]\simeq J[2]^\vee$
is its own dual because of the Weil pairing (we recall that for a vector space $V$ the automorphism $\sigma\in G_\Q$ acts on the dual $\varphi \in V^\vee$ of $V$ by $\sigma \cdot\varphi := \varphi \circ \sigma^{-1}$), we obtain the following short exact sequence of Galois-modules
$$ 0 \ra J[2] \ra E^\vee \ra R^\vee \ra 0.$$
Taking Galois cohomology for $k=\Q$ or $k=\Q_p$ now leads to 
$$ 0 \rightarrow J[2](k) \rightarrow E^\vee (k) \rightarrow R^\vee (k) \ra \HJk \overset{\alpha_k}{\rightarrow} \Hop^1(k,\E) \ra \Hop^1(k,\Rw)$$
and therefore to the exact sequence
$$ 0 \ra R^\vee(k) /E^\vee(k) \ra \HJk \overset{\alpha_k}{\ra} \Hop^1(k,\E).$$
We call the last arrow in the present article the \emph{global} or \emph{local cohomology map} $\alpha_k$, according to whether $k=\Q$ or $k=\Q_p$. The full cohomology group $\Hop^1(k,E^\vee)$ is not explicit itself but a sufficient large subgroup is as given in the following property. We define $L=\operatorname{Hom}_{G_\Q-sets}(\Delta,\Qb)$ to be the \'etale algebra associated to $\Delta$. We denote $L_k:=L\otimes_\Q k$ the $G_{\Q_p}$-\'etale algebra associated to $\Delta$ and $Br(k):= \Hop^2(k,k^*)$ the Brauer group of $k$ and $\mu_2^\Delta:=(\F_2^\Delta)^\vee$ the dual of $\F_2^\Delta$.

\begin{prop} The Galois-modules $\E$ and $\mu_2^\Delta/\mu_2$ are isomorphic.\\
    We have the isomorphism $\etalek \simeq \ker(\Hop^1(k,\mu_2^\Delta/\mu_2) \ra Br(k))$. Moreover, it contains the image of $\HJk$ by $\alpha_k$.
\end{prop}

\begin{proof} 
\cite[section 6.3.3. and Prop 6.10.]{BPS}.
\end{proof}

In case of a True Descent, the situation is quite similar. 

\begin{prop} 
    We have the isomorphism $L_k^*/L_k^{*,2} \simeq \Hop^1(k,\mu_2^\Delta)$.
\end{prop}

Note that $\etalek$ is just the space $L_k^*/L_k^{*,2}$ quotiented by the diagonal action of $k^*$.\par
We will later see in another section how to obtain $L$ as an explicit quotient algebra $\Q[y]/\chi(y)$. So finally, we have the following exact sequence
$$ 0 \ra R^\vee(k) /E^\vee(k) \ra \HJk \overset{\alpha_k}{\ra} \etalek \subset \HDeltak,$$
where we know how to make $\etale$ explicit in practice. As we shall soon see, it is the ambient space in which the \emph{Fake Selmer group} lives. The most visible difference between Fake and True Descent is that this ambient space becomes $L^*/L^{*,2}$ with the appropriate \'etale algebra $L$ like in the elliptic case we have exposed.

\subsection{Computing the Fake Selmer Group $\SelF$}

We note $\res_p$ for the natural maps $\HJ \ra \Hop^1(\Q_p,J[2])$ induced by the restriction to a decomposition group, $\gamma_p: \JQp \hra \Hop^1(\Q_p,J[2])$ the connecting map of the cohomology sequence stemming from $0 \ra J[2] \ra J \ra J \ra 0$, and $\alpha_k : \HJk \ra \HDeltak$ the map defined above. We also note $C_k := \alpha_k \circ \gamma_k$ and call it \emph{the global} (respectively \emph{local at $p$}) \emph{descent map}. This map is of crucial importance for explicit 2-descent. The main reason for the full subtilty of the Fake Descent setup's definition (See \cite[Def 6.7.]{BPS}) over the framework we set so far is that it ensures we can break down $C_k$ to an explicit computable map. In other words, it ensures that \cite[Prop 6.10.]{BPS} is satisfied. We will see how to compute the explicit components of the map $C_k$ in the following chapter. As the image of $C_k$ belongs to the subgroup identified with $\etalek$ inside $\HDeltak$, we will often not distinguish the two maps $\Jk \ra \HDeltak$ and $\Jk \ra \etalek$. We warn the reader that we will sometimes also identify $C_k$ with its composition with the canonical projection $J(k) \ra \Jk$. For what's following, it is helpful to keep the following commutative diagram in mind.
\[
\begin{tikzcd}
    J(\Q_p)/2J(\Q_p) \arrow[d, hook, "\gamma_p"] \arrow[r, "C_p"] & L_p^*/L_p^{*,2} \arrow[d, hook]  \\
    \Hop^1(\Q_p,J[2]) \arrow[r, "\alpha_p"] & \Hop^1(\Q_p,\mT)
\end{tikzcd}
\]

\subsubsection{The Selmer group $\Sel$}

We recall that the 2-Selmer group $\Sel$ of the Jacobian is defined by
$$ \Sel := \{ \zeta \in \HJ ~ | ~ \forall p ~ \text{place of $\Q$, } \res_p(\zeta) \in \gamma_p (\JQp) \}.$$

\subsubsection{The Fake Selmer group $\SelF$}

Now that $\alpha_k$ allows us to make $\HJk$ explicit, it would be natural to try to compute the image $\alpha_\Q(\Sel) \subset \etale\subset \HDelta$. However, it is actually easier to compute the subgroup of $\etale$ cut out by the same conditions, and which is called the \textit{Fake Selmer group} by \cite{BPS}:
$$ \SelF:= \{ \zeta \in \etale ~ | ~ \forall p \text{ place of $\Q$, }\res_p(\zeta) \in C_p(\JQp) \}. $$
We have the inclusion $\alpha_\Q(\Sel) \subset \SelF$.

\subsubsection{The approximation $S_\Tset$}

The ambient spaces $\etalek$ are explicit in that we can represent their elements effectively, but still of infinite dimension. We therefore consider a smaller space to start cutting out local conditions in. Let $\Sset $ be the set of places of bad reduction for the Jacobian along with the places above $2$ and the archimedean ones. As we now show, the elements of $\SelF$ are unramified (See Definition \ref{def unramified}) at least at every prime of good reduction. See \cite[Lem 7.1.]{BPS} for a generalization of the following fact.
\begin{prop}\label{Selunramif}
    Let $p\neq 2$ be a prime of good reduction.\\
    Then, $\gamma_p(J(\Q_p)/2J(\Q_p)) = \Hop^1_{unr}(\Q_p,J[2]).$
\end{prop}
\begin{proof}
    We first show the inclusion $\subset$. By the commutativity of the underlaying diagram, the restriction of a cocycle $\xi \in \gamma_p(J(\Q_p)/2J(\Q_p))$ to the Inertia group $I_p$ belongs to $\gamma_p(J(I_p)/2J(I_p))$. The latter group is trivial since, for a prime of good reduction, the action of $I_p$ over $J$ is. We recall that $J(I_p)/2J(I_p) =J(\Qpb)/2J(\Qpb)=0$ since $\Qpb$ is algebraic closed.\\
    From the restriction-inflation sequence, we have that $\Hop^1_{unr}(\Q_p,J[2])=\Hop^1(<\Frob_p>,J[2]^{I_p})$. Because of Proposition \ref{LOTC ramification}, the action of $I_p$ over $J[2]$ is trivial. It implies that $\Hop^1(<\Frob_p>,J[2]^{I_p})=\Hop^1(<\Frob_p>,J[2])$. We conclude by cardinality. For a prime $p\neq 2$, because $<\Frob_p>$ is cyclic, we have $|\Hop^1(<\Frob_p>,J[2])|=|J[2]^{\Frob_p}|=|J[2](\Q_p)| = |\gamma_p(J(\Q_p)/2J(\Q_p))|$ (See the Proposition \ref{local card} below for the last equality).
\end{proof}
Thus, $\SelF$ belong to the $\F_2$-vector space $\sLS$ of elements of $\etale$ unramified everywhere outside $\Sset$, which is \emph{finite-dimensional} (See \cite{BPS} section 7.3.). Moreover, every local condition can, in theory, be explicitly checked for an element $\zeta \in \sLS$, as doing so only requires to verify that $\zeta$ lies in $J(\Q_p)/2J(\Q_p)$ which is a $\F_2$-vector space of known finite dimension. 
\begin{prop}\label{local card}
    Let $k_p$ be the completion of a global field $k$ at a prime $p$.\\
    Then $\dim_{\F_2} J(k_p)/2J(k_p) = \dim_{\F_2} J[2](k_p) + g\delta_{p=2}$.
\end{prop}
\begin{proof}
    See \cite[Prop 12.10.]{PoonenSchaefer}.
\end{proof}
This is why, for any set of places $\Tset$, following BPS we note
$$ S_\Tset := \{\zeta \in \sLS ~ | ~ \forall p \in \Tset, ~ \res_p(\zeta) \in C_p (\JQp )\}$$
the subgroup of $\etale$ with local conditions at $\Tset$ and unramified everywhere outside $\mathcal{S}$. Finally, for any \textit{finite} $\Tset$, the set $S_\Tset$ is completely computable (in theory) as a subgroup of an explicit vector space of finite dimension defined with a finite number of local checkable conditions. Rather than $\SelF$ itself, the group $S_\Tset$ is the actual group that we aim and reach computationally.\bs

Nevertheless, we always have the inclusion $\SelF \subset S_\Tset$. Moreover, we observe that very few places are needed to achieve $\SelF = S_\Tset $ in practice. The reader should notice that in contrast with the (hyper)elliptic case where $\Hop^1(\Q,\Rw)\simeq \Q^*/\Q^{*,2}$, we cannot distinguish the elements of $\sLS$ coming from $\HJ$ because the map $\etale \subset \Hmu \ra \Hop^1(\Q,\Rw)$ is not explicit. As for $\HJ$, we do not even know how to represent the elements of $\Hop^1(\Q,\Rw)$ on a computer. Due to this defect, we may observe $\SelF \neq S_\Sset$ for Jacobians of general curves. This means that local conditions at good places may matter as well. Indeed, the local condition at a place $p$ of good reduction does not match with being unramified at $p$ in general. It imposes a bigger constraint and we cannot bridge the gap. Formally, denoting $I_p \subset G_{\Q_p}$ the inertia group at $p$ of good reduction, from Proposition \ref{Selunramif} we have 
$$C_p(J(\Q_p)/2J(\Q_p)) \subset \ker(\Hop^1(\Q_p,\mT) \ra \Hop^1(I_p,\mT))$$
with no necessary equality. However 
$$C_p(J(\Q_p)/2J(\Q_p)) = \ker(\Hop^1_{unr}(\Q_p,\mT)\ra\Hop^1(\Q_p,\Rw))$$
but the map is not computable in general. We can still be ensured to recover $\SelF$ with a finite number of local conditions if we take $\Tset$ to contain $\Sset$ and enough places of good reduction so that the corresponding Frobenius elements reach all of the conjugacy classes of the Galois group of $J[2]$ (See \cite[Thm 10.9]{BPS} for a formal statement). Unfortunately, it usually makes $\Tset$ huge in a computational point of view and unpracticable. In practice, we start with a few places and make $\Tset$ bigger until we are satisfied. It is not rare that $\Tset = \{2\}$ is enough to meet the lower bound on $dim_{\F_2} \JQ$.\bs

In opposite to the elliptic case as well, the short exact sequence $0 \ra J[2] \ra \E \ra \Rw \ra 0$ does not split. A consequence is that the maps $\alpha_k : \HJk \ra \HDeltak$ and thus the maps $C_k : J(k)/2J(k) \ra L_k^*/L_k^{*,2}k^*$ are not injective in general either globally for $k=\Q$ or locally over the $p$-adic for $k=\Q_p$. We therefore need to control those kernels. This is the topic of the next section, following \cite{BPS}.

\subsection{Controlling cohomology kernels $\kF(k)$ and descent kernels $\ker(C_k)$}\label{sectionKernel}

See \cite[section 10.]{BPS}. As we said, the maps $C_k : \Jk \ra \etalek$ are not injective in general for $k=\Q$ or $k=\Q_p$. This leads to two distinct problems.

\subsubsection*{Globally}
The cardinal of $\SelF$ is not a bound of the cardinal of $\Sel$. To obtain such a bound, we have to control the kernel of the map $\alpha_\Q$ not only over $\HJ$ but over the subgroup $\Sel$ for a better result. This would provide the requested inequality $|\Sel| \leq |\SelF| \times|\ker((\alpha_{\Q})_{|\Sel})|$.\par
The reader should note that it is not always an equality as $\alpha_\Q(\Sel) \subset \SelF$ with no reciprocal inclusion guaranteed. Indeed, as we already said, in contrast to the (hyper)elliptic case where $\Rw = \mu_2$, the last map of the exact sequence $\HJ \ra \Hmu \ra \HR$ is not explicit anymore. Therefore, we cannot exactly distinguish the elements in $\SelF$ which are not in $\alpha_\Q(\HJ)$. See the discussion above. 

\subsubsection*{Locally}
If it is not that complicated to produce points in the local image $C_p(\JQp)$, it is however more difficult to know when we have reached the full image. We know that the $\F_2$-dimension of $J(\Q_p)/2J(\Q_p)$ is $\dim_{\F_2}(J[2](\Q_p)) + g\delta_{p=2}$ but this does not tell us the dimension of the local image when $C_p$ is not injective.\bs

We now discuss these concerns in more detail. Recall that $\E \simeq \mT$. We assume that we have computed $\SelT$ for some set of places $\mathcal{T}$ such that we are rather convinced that it is exactly $\SelF$. 
We now describe how to recover a bound of $|\Sel|$ using the methods of \cite{BPS}. 
We refer to section \ref{section explicit Galois} for explicit computations. Following the notations of BPS, the following sequence is exact

\begin{tikzcd}
    &&&& J(k)/2J(k) \arrow[d, hook, "\gamma_k"] \arrow[rd, "C_k"] & & \\
    0  \arrow[r] & J[2](k) \arrow[r] & \E(k) \arrow[r] & \Rw(k) \arrow[r] & \HJk \arrow[r, "\alpha_k"] & \Hop^1(k,\E) ,
\end{tikzcd}

\bigskip

and we have the chain of inclusions of global spaces:
$$ \alpha_\Q (\Sel) \subset \wSel  \subset \SelF \subset \SelT \subset \sLS \subset \Hop^1(\Q,\E)$$
with $\wSel := \SelF \cap \alpha_\Q (\Hop^1(\Q,J[2])$. In other terms, we want to control the gap between $\Sel$ and its image $\alpha_\Q (\Sel)$ and the gap between the latter and $\SelF$.\\

We set in the present article the following definitions.

\begin{de}\label{def CohoKernel}
According to $k=\Q_p$ or $k=\Q$, we call the local or global \underline{cohomology kernel}:
$$ \kF (k) := \ker ( \alpha_k) \simeq \Rw (k) / \E (k).$$
It should not be confused with the \underline{global descent kernel}:
$$ \sKF := \kF(\Q) \cap \Sel = \ker(\alpha_{|\Sel})$$
and the following \underline{local descent kernels}. Locally, for any prime $p$:
$$ \ker (C_p) = \kF(\Q_p) \cap  \gamma_p(J(\Q_p)/2J(\Q_p)).$$
\end{de}

Following BPS, we set
$$ W_p := \kF(\Q_p) / \ker (C_p).$$

\subsubsection{The gap $|\Sel| - \alpha_\Q (\Sel)$: Computing the global descent kernel $\sKF$}

\begin{prop}\label{Wp}
    For any prime $p\neq 2$ of good reduction, $W_p = 0.$
\end{prop}

\begin{proof} (See also \cite{BPS} Corollary 10.6)\\
    From Proposition \ref{Selunramif}, it is enough to show that $\ker(\alpha_p)$ is unramified. Now, the restriction of a cocycle $\xi \in \ker(\alpha_p) = \kF(\Q_p)$ to the Inertia group $I_p$ belongs to $\Rw(I_p)/\E(I_p)$. But the latter group is trivial since the action of $I_p$ over $\OTC$ is (See Proposition \ref{LOTC ramification}).
\end{proof}

In particular, the space $\prod_{\forall p} W_p = \prod_{p\in \Sset} W_p$ is finite. Let $ \kappa : \kF(\Q) \ra \prod_p W_p$ be the natural projection map.

\begin{prop}
    The following sequence is exact:
    $$ 0 \ra \ker (\kappa) \ra \Sel \overset{\alpha_\Q}{\ra} \wSel \ra \coker (\kappa).$$
\end{prop}

\begin{proof}
    (Snake lemma). See \cite[Thm 10.14]{BPS}.
\end{proof}

This has two explicit consequences. First, $\sKF \simeq \ker(\kappa) \text{ is computable in practice.}$ Secondly, the gap $\alpha_\Q(\Sel) - \wSel$ can sometimes be detected when trivial.

\begin{prop}
    If $\coker( \kappa) = 0$, then $\alpha_\Q(\Sel) = \wSel.$
\end{prop}

\subsubsection{Computing $W_p$ and local images $\im(C_p)$}\label{section computation local image}

Assume that we have found $n$ different points in $\im(C_p) = C_p( \JQp)$ for a prime $p$. We present methods to check if we have found the whole image. We then deduce the space $W_p$.
\begin{itemize}
    \item If $\kF(\Q_p) = 0$, then $C_p$ is injective (because $\alpha_p$ is). It is enough to check that $n= |J(\Q_p)/2J(\Q_p)| = |J[2](\Q_p)| \times 2^{g\delta_{p=2}}$ to know when we are done. It follows that $W_p = \kF (\Q_p)$.
    \item Fortunately, it can still happen that $\alpha_p$ is not injective but $C_p$ is.
    \item In the (not so frequent) general case, we can compute (with some time) $\im(C_p)$ using the "residue discs method" of BPS (See \cite[section 11.1]{BPS} or the section \ref{section reducing kernel} of the present article) taking sufficiently enough points in all extensions of $\Q_p$ of bounded degree. It directly gives $\ker(C_p)$ and we can therefore compute $W_p$. It is generally less and less practicable as the genus increases.
\end{itemize}

\subsubsection{The gap $\wSel - \SelF$}

In contrast with the (hyper)elliptic case, we cannot easily tell which elements of $\Hop^1(k,\E)$ come from $\Hop^1(k,J[2])$. As we said, this is mainly because the map $\HE \ra \HR$ is not explicit.

\begin{de}
    For any Galois-module $M$, we define the Shafarevich group-like space\\
$\Sha^1(k,M):= \bigcap_v \ker \big( \Hop^1(k,M) \ra  \Hop^1 (k_v,M) \big)$.
\end{de}

\begin{prop}
    The following sequence is exact
    $$ 0 \ra \wSel \ra \SelF \ra \Sha^1(\Q,\Rw).$$
\end{prop}

\begin{proof}
    See \cite[Lem 10.17.]{BPS}.
\end{proof}

As often with Shafarevich groups, we expect it to be frequently trivial. It is indeed often the case in genus $g=3$ (See the discussion \cite[section 12.7.(ii)]{BPS}).



\subsubsection{Summary}

\begin{itemize}
    \item $|\Sel | = |\alpha_\Q(\Sel)| \times |\sKF|$, with $\sKF$ explicit.
    \item $|\alpha_\Q (\Sel) | \leq |\wSel|$, with equality if $\coker(\kappa)=0$ (explicit).
    \item $|\wSel| \leq |\SelF |$, with equality if $\Sha^1(\Q,\Rw) =0$. 
    \item $|\SelF| \leq |\SelT |$, with equality if $\mathcal{T}$ contains enough places.
\end{itemize}
In particular, we obtain the following computable bound for the Selmer group:
$$ |\Sel| \leq |\SelT| \times |\sKF|$$
with equality in particular if $\mathcal{T}$ contains enough places, $\coker(\kappa)=0$ and $\Sha^1(\Q,\Rw) =0$. The condition $\coker(\kappa)=0$ is quite easy to guarantee when it happens but the author has chosen not to check the two others.


\section{Performing generalized explicit 2-descent in practice}

In this chapter we describe an automatized algorithm implemented by the author designed to perform generalized explicit descent in the BPS meaning (See \cite{BPS}) on the Jacobian $J$ of a given "nice" (See Notations in the introduction) curve $X$ defined over $\Q$ and equipped with a rational point. The algorithm requires the following inputs:
\begin{enumerate}
    \item An equation $f(x,y)\in \Q[x,y]$ whose associated (singular) projective model in $\mathbb{P}^2$ is birational to the curve $X$ we are interested in.
    \item A rational point $P_0$ on the model $f$ (or a degree 1 rational divisor).
    \item A description of the Galois group $G$ of $J[2]$ as embedded in $\Sp$.
\end{enumerate}
We can try to search input 2. by brute-force if not provided. Methods to obtain input 3. from the data we compute anyway are described in \cite{BPS}. The algorithm first acquires the OTC and their Galois action and then performs generalized 2-descent as described earlier. If successful, it returns 
\begin{enumerate}
    \item $B$ a bound for $|\Sel_\Q(J)|$
    \item $B-\dim_{\F_2}(J[2](\Q))$ a bound for $\rank_\Q(J)$
\end{enumerate}

It heavily uses another algorithm \cite{TitiOTC} implemented by the author which computes the (Odd) Theta Characteristics of the curve and their attached invariants. Both algorithms are implemented in the free and open source software language PARI/GP \cite{PARI2}.\bs

We discuss our implementation of 2-descent in this section. We first gather some facts on the prerequisites for 2-descent and how we compute them. 
We then give some insights on the technicalities upon the implementation we made.

\subsection{Computing Theta Characteristics}

What usually prevents from starting a 2-descent on the Jacobian of a non-hyperelliptic curve is the difficulty to compute the necessary OTC/2-torsion invariants on which the various maps and spaces described earlier are based on. Generalized 2-descent in the BPS meaning requires at least the following.
\begin{enumerate}
    \item The \'etale algebra $L$ attached to $\OTC$.
    \item The local and global descent maps $C_k : J(k) \ra J(k)/2J(k) \ra \Hop^1(k,J[2]) \overset{\alpha_k}{\ra} \Hop^1(k,\mu_2^{\OTC}/\mu_2)$ for $k=\Q$ and $k=\Q_p$ for all $p\in\Tset$.
\end{enumerate}

In the (hyper)elliptic case $E : y^2=f(x)$, the \'etale algebra $L$ is easily obtained as $\Q[y]/f(x)$ and the map of components $(x-\gamma_i)_i$ suits as descent map if $\gamma_i$ are the roots of $f$. Shortly, it comes from the fact that the principal divisors $(\divi(x-\gamma_i))_i$ express the divisors of the OTC (or more precisely of a OTC Galois-subset spanning the whole 2-torsion). This is the reason why $p$-descent is usually performed on $p$-covers of $\mathbb{P}^1$ (See \cite{PoonenSchaefer}). These invariants are more complicated to compute in general. Given the inputs mentioned earlier, the algorithm \cite{TitiOTC} of the author returns by design the mentioned data required to carry out 2-descent. We present this below. But first we make a remark on the input of a rational point.

\subsubsection{Hypothesis 10.1}

The input of a rational point $P_0$ (or a degree 1 rational divisor) is a necessary practical assumption so that the algorithm \cite{TitiOTC} computing the OTC runs. This is usually considered as a strong assumption as it is complicated to check the existence of such a point in general. However, it ensures that the Hypothesis 10.1 required in \cite{BPS} and the identification $J \overset{\Q}{\simeq} \Pic_X^0$ over $\Q$, essential for Galois to wellbehave, are satisfied.

\subsubsection{Computing the required data for 2-descent}

In \cite{BPS}, the authors essentially derive the required data $1.$ and $2.$ from the bitangents of the curve in its canonical model (See lemma 12.6.). When $g=3$, a smooth plane quartic is naturally in its canonical model and the equations are small enough so that Groebner basis computations terminate. Stoll and Fu conducted the same process on a genus $4$ curve in \cite{StollFu}. Although ETC can be expressed by OTC, this technique \textit{a priori} cannot extend directly to Even Theta Characteristics as they have no geometrical counterpart.\bs

The algorithm \cite{TitiOTC} of the author works differently. If it takes advantage of the geometrical interpretation of the OTC by the bitangents, it does not rely on it and therefore applies equally to OTC and ETC.\par
The algorithm \cite{TitiOTC} computes the image of an injective Galois-set morphism $\pi : \OTC \ra \Qb$ and the associated polynomial $\chi_\pi(y):=\Pi_{\beta \in \OTC} (y- \pi(\beta))$. The polynomial $\chi_\pi$ is consequently in $\Q[y]$ and induces the \'etale algebra $L:=\text{Hom}_{G_\Q-set}(\OTC,\Qb)\simeq \Q[y]/\chi_\pi(y)$.\par 
The algorithm \cite{TitiOTC} also explicitly computes the components of the descent maps $C_k$. For each Galois orbit representative $\beta$ of $\OTC$, we compute a function $f_\beta$ on the curve called \textit{Theta Fraction} in the article. Its principal divisor is
$$ (f_\beta) = 2\underline{\beta} - 2K_X$$
where $\underline{\beta}$ is a divisor in the class of $\beta$ and $K_X$ is fixed in the canonical class. The choices of the fixed canonical divisor $K_X$ and of the divisors $\underline{\beta}$ have no theoretical incidence. For computational efficiency, we choose $\underline{\beta}$ to be the unique \textit{effective} divisor in the class $\beta$ when possible, i.e. when the Riemann-Roch dimension $h^0(\beta)$ equals 1. As suggested by \cite[Thm 2.17]{ModuliSpaceEffective}, this is often the case. 
According to \cite[Prop 6.10.]{BPS}, the map whose components are these Theta Fractions is suitable for the global descent maps $C_\Q$ under the canonical projection $L \ra \etale$. We just compose with the $p$-adic embeddings to obtain the local descent maps $C_p$.\bs

It can be convenient to keep in mind the following commutative diagram with exact rows. For $k=\Q,\Q_p$

\begin{tikzcd}
    && J(k)/2J(k) \arrow[d, hook, "\gamma_k"] \arrow[r, "C_k"] & L_k^*/L_k^{*,2}k^{*,2} \arrow[d, hook] & \\
    0  \arrow[r] & \Rw(k)/\E(k) \arrow[r] & \Hop^1(k,J[2]) \arrow[r, "\alpha_k"] & \Hop^1(k,\mT).
\end{tikzcd}

\subsubsection{Outputs of the TC algorithm}

The algorithm \cite{TitiOTC} of the author computes the required \'etale algebra $L$ as well as the local and global descent maps $C_k$.\bs


It also returns the following useful data. It will ease the further coming 2-descent.
\begin{enumerate}
    \item The \'etale algebra $L_{tors}$ of $J[2]$.\\
    By the same process mentioned earlier, we compute a polynomial $\chi_{J[2]}$ realizing $L_{J[2]}$ as a quotient algebra of $\Q[y]$. It allows us to know easily the dimensions of the local spaces $J[2](\Q_p)$ framing the ones of the local images as well as the global one $J[2](\Q)$ artificially increasing the dimension of the Fake Selmer group.
    \item The genus $g$ of the curve.
    \item The set of bad prime of the curve.\\
    These are useful to determine as few relevant local conditions as possible.
    \item The Quadratic Theta Forms of the OTC.\\
    We use this knowledge to compute explicitly the Galois action on the local and global cohomology kernels $\kF$.
    \item The Steiner Complexes attached to the OTC.
    \item The \'etale algebra $L_{TC}$ of the whole set $\TC$ of the Theta Characteristics.
    \item The Frobenius over $J[2](\Qpb)$ for any chosen prime $p$ of good reduction.
\end{enumerate}

Not returning these outputs does not really spare time as their computation is either essential or heavily based on the precomputions the main outputs demand.

\subsection{Implementing 2-descent}

We implemented 2-descent in the free and open source software language PARI/GP \cite{PARI2} \cite{PARI2}.

\subsubsection{Possible failures}

We quickly describe what may fail in a $2$-descent.\par
First, we may fail to obtain the OTC data or the Galois group, which are required to perform $2$-descent in the way we presented it. For our computations, a rational point on the curve is required. However, the author observed a high success for the curves he ran the algorithms on. Otherwise, see the techniques presented in $\cite{BPS}$.\par
The $2$-descent itself is a method aiming at bounding the rank. In particular, it is not guaranteed to reach exactly the rank or give a satisfactory bound. It may due to a non-trivial Tate-Shafarevich group $\Sha^1(\Q,J)[2]$. In that case, the Selmer group $\Sel$ does not match $J(\Q)/2J(\Q)$. Even if it does, the cardinal $|\SelF|\times |\sKF|$ may fail to approximate correctly the one of $\Sel$. Either because we were not able to compute completely $\sKF$ but an upper bound or because the map $\alpha_\Q : \Sel \ra \SelF$ is not surjective. In other words, looking at $\Hop^1(\Q,\E)$ rather than $\Hop^1(\Q,J[2])$ may introduced new elements that we were not able to distinguish from the ones coming from $\Sel$. However, as showed in the Examples section, it works quite well in practice.\par
We may also fail to compute the class groups and units groups of the fields composing the \'etale algebra $L$. We dive more into this condition in the next section.

\subsubsection{The need for splitting orbits}

We use PARI/GP \cite{PARI2} algorithms to compute the classgroups and basis of the unit groups of the fields of the \'etale algebra.\par
As we can't afford to compute them for fields of degree $|\OTC|=120$ when $g=4$ for instance, we need the orbits to split in order to perform 2-descent if $g\geq 4$. Along with having a rational point, they are the main conditions on the curve required to perform 2-descent in the way we implemented it. We give here evidences for orbits to split. Namely, we prove that non-trivial rational $2$-torsion, automorphisms of the curve or \'etale morphisms defined over the groundfield $k$ split the TC orbits. We do not need to be aware of any of the following situations to benefit from it.

\begin{prop}\label{TCsplit}
    Let $X$ be a "nice" curve defined over $k$ and $J$ its Jacobian (with $g>1$).\\
    If $J[2](k)$ is non-trivial, then the Galois action of $G_k$ is not transitive over $\OTC$ nor $\ETC$.\\
    Moreover, in that case, the biggest orbit is respectively of size at most $2^{2g-2}$ and at most $2^{2g-2} + 2^{g-1}$.
\end{prop}

\begin{proof}
There exists a non-trivial rational 2-torsion point $P_0$ by assumption. For any $\beta\in \TC$, the TC $\beta$ and $\beta+P_0$ have the same orbit length because they share the same stabilizer. We can choose $\beta$ such that $\beta$ and $\beta + P_0$ do not have the same parity. Indeed this is the case if and only if $q_\beta(P_0)=1$. But, up to fixing a QTF $q_0$, every quadratic theta form can be expressed as $q_0 +\<P,\cdot\ran_W$ for a certain $P\in J[2]$. Indeed, one can check that this gives quadratic forms. By cardinality, we reach them all. Conversely, every $P\in J[2]$ defines a quadratic theta form this way. As $P_0$ is non-trivial, we can find $P\in J[2]$ such that $\<P,P_0\ran_W=\qs(P_0) +1 \in \F_2$. The TC $\beta$ corresponding to the QTF $q_0 + \<P,\cdot\ran_W$ has the desired property.\par
Now, since $\OTC$ and $\ETC$ do not have the same cardinal and that both sets are preserved by the Galois action, the latter can't be transitive over $\ETC$.\par
We have to be more careful to show that the orbits of $\OTC$ split as well. There are exactly $2^{2g-1}$ points $P\in J[2]$ such that $\<P,P_0\ran_W=q_0(P_0) +1$. Each of them defines a pair of TC $(\beta,\beta +P_0)$ spanning one orbit in $\OTC$ and one in $\ETC$ of same length. We observe that $q_{\beta+P_0}(P_0) = q_\beta(P_0)=1$. So there is a set $\Delta_c$ of cardinal $2^{2g-2}$ formed by such pairs which is stable under Galois and partitioned in orbits of lengths that are included in the set of lengths of orbits of $\OTC$ (and of $\ETC$). Since $|\OTC|=2^{g-1}(2^g-1)>2^{2g-2}$, the OTC split in two non-trivial Galois-sets according to the parity of their sum with $P_0$. For the record, one of the Galois-set is given by the OTC of the pairs in $\Delta_c$ and the other is exactly formed by the OTC belonging to the Steiner Complex of $P_0$. Moreover, the biggest of this two Galois-sets is of size $max(|\OTC|-|\Delta_c|,|\Delta_c|)=2^{2g-2}$. This gives a bound on the size of the biggest OTC orbit. Similarly, for the ETC, the biggest orbit is of size at most $max(|\ETC|-|\Delta_c|,|\Delta_c|)=2^{2g-2}+2^{g-1}$.
\end{proof}

The last property is true over $\Q$ but over $\Q_p$ as well, where there is more often non-trivial rational torsion. This is good news as we could be willing to compute the Galois groups of the fields of the local \'etale algebras. For that purpose, the smaller the size of the orbits are, the better.\par
\begin{example}
    From the algorithm \cite{TitiOTC}, the smooth plane quartic defined by $ y^3=x^4+6x^2+1$ have orbits decomposition $ \OTC: 1_2 ~ 2_1 ~ 8_1 ~ 16_1$ and $J[2]: 1_1 ~ 2_1 ~ 4_1 ~ 8_3 ~ 16_2$. The bounds given above therefore seem to be sharp.
\end{example}
\begin{rem}
    For $g=4$, the given bounds for the biggest orbit are respectively $64$ and $72$, which would make the computation of the classgroup and units group of the corresponding field still quite involved. If we are lucky, the orbits split more than this.\bs
\end{rem}

On the other hand, curves with non-trivial rational automorphisms have better chance to have split $\OTC$-orbits as well. Indeed, given a group $G\subset \Sp$, we let $\operatorname{Com}_{\operatorname{Sp}} (G):=\{\sigma\in \Sp ~|~\forall g\in G, ~\sigma\circ g =g\circ \sigma \}$ be its centraliser. Let $\Aut_k(X)$ be the group of automorphisms of the curve defined over the base field $k$. It extends to $J$, commutes with the group law and therefore preserves the $2$-torsion. Let $\Aut_k(J[2])$ be the group of restrictions to $J[2]$ of endomorphisms of the Jacobian defined over $k$. Then, we have 
$$\im(\rpztors) \subset \operatorname{Com}_{\operatorname{Sp}}(\Aut_k(J[2])).$$ 
In particular, if $\Aut_k(X)\neq \varnothing$, it is more likely that $\im(\rpztors)$ be strictly smaller than the expected group $\Sp$, and therefore the TC orbits are more likely to split.\bs

Here is another evidence for splitting orbits.
\begin{prop}\label{etaleTC}
    Let $X$ be a "nice" curve defined over $k$ of genus $g$. Assume that $X$ admits an \'etale morphism $\pi: X \ra \widetilde{X}$ defined over $k$ of degree $d$ toward a curve $\widetilde{X}$ of genus $\widetilde{g}$.\\
    If $0<\widetilde{g}<g$ and $d$ is odd, then the orbits of $\TC(X)$ split.\\
    Moreover, if there is no non-vanishing ETC, the orbits of $\OTC(X)$ split as well.  
\end{prop}

\begin{proof}
    The map $\pi$ induces a \textit{norm map} $\pi_* : \Pic_X \ra \Pic_{\widetilde{X}}$ as the push forward of $\pi$. We also consider the pull back $\pi^*: \Pic_{\widetilde{X}} \ra \Pic_X$ of $\pi$. The two maps satisfies $\pi_* \circ \pi^* = [d]$. 
    Since $\pi$ is defined over $k$, so is $\pi^*$. Moreover, since it is \'etale, $\pi^*$ also sends the canonical class of $\widetilde{X}$ over the one of $X$ and we therefore obtain a restrained map $\pi^* : \Delta_{TC,\widetilde{X}}  \ra \Delta_{TC,X}$ defined over $k$ and Galois equivariant.\par
    Assume there are two TC $\beta, \beta'$ of $\widetilde{X}$ that have the same image under $\pi^*$. Then, composing with $\pi_*$, we observe that $[d](\beta - \beta') =0$. But $\beta - \beta'$ is a 2-torsion point and $d$ is odd. So we have $\beta = \beta'$.\par
    Since $\pi^*$ is injective on the TC, the image of $\Delta_{TC,\widetilde{X}}$ gives strictly non-trivial suborbits of $\Delta_{TC,X}$ because $0<\widetilde{g}<g$.\par
    In the usual case where there is no non-vanishing ETC on $X$, the OTC are the only TC to have non-trivial Riemann-Roch spaces. Since for any TC $\widetilde{\beta}$ of $\widetilde{X}$ the map $\pi^*$ gives rise to an injective map $\Hop^0(\widetilde{\beta})\hra \Hop^0(\pi^*(\widetilde{\beta}))$, we have that $\Delta_{Odd,\widetilde{X}}$ is sent into $\Delta_{Odd,X}$. We therefore further deduce that $\Delta_{Odd,X}$ splits. Galois may remain transitive over $\Delta_{Even,X}$ if $\pi^*$ sends the ETC of $\widetilde{X}$ into the OTC of $X$.
\end{proof}

Note that the above property cannot apply when $X$ is an elliptic curve, as the Riemann-Hurtwitz formula would impose $g= \widetilde{g}$.

\subsubsection{Computing local images $\im(C_p)$}

For computing the local image $C_p(J(\Q_p))$, we produce points $P$ on $X(\Q_p)$ with finite $p$-adic approximations and consider differences $D:=P-P'$. We then approach the coefficients of a representative of $C_p(D)$ in $L_p$ with rational numbers in $\Q$ and obtain the image in $L_p^*/L_p^{*,2}$ via restriction from $L_\Q$ by the appropriate primes.\par
Related to that last step, we recall that for any extension $K/\Q_p$, the quotient $K^*/K^{*,2}\simeq \Z/2\Z \times O_K^*/O_K^{*,2} \simeq \Z/2\Z \times \F_K^*/\F_K^{*,2} \times (1+(\pi))/(1+(\pi))^2$ is a finite $\F_2$-vector space and that the projection $K^*\ra K^*/K^{*,2}$ of an explicit element can be computed from its valuation and local restriction in $\F_K$. So do elements in $L_p^*/L_p^{*,2}$. Computing the projection $K^*\ra K^*/K^{*,2}$ is quite elementary for $p\neq2$ where $(1+(\pi))/(1+(\pi))^2 \simeq 1 $ is trivial. For $p=2$, where $\F_K^*/\F_K^{*,2}\simeq1$, we use Stoll's algorithm \cite[below Lem 3.1]{StollCH}.\par
See section \ref{section computation local image} for details on the (classic) method we use.

\subsubsection{Computing $\SelT$}\label{section ST}

We compute $\SelT$ by cutting out the local conditions starting from $\sLS \simeq (\etale)_S \hookrightarrow H_S^1(\Q,\mT)$ as in \cite[section 7.3]{BPS}. In practice we compute $\sLS$ inside $\LS \simeq H^1_S(\Q,\muT)$ and recover the latter from $L(2,\Sset')$ with $\Sset' \supset \Sset$ large enough to contain primes that span the classgroups of $L$. In that case, the vector space $L(2,\Sset') \simeq \Oc_{L,\Sset'}^*/\Oc_{L,\Sset'}^{*,2}$ is easier to represent. See \cite{PoonenSchaefer} section 12. for computing basis of these different spaces efficiently.\bs

Ideally we would have asked the elements of the Fake Selmer group to be also in\\ $\im(\alpha_\Q) = \ker(\Hmu\ra \Hop^1(\Q,\Rw))$ but as we already said this cohomology map is not explicit. So, according to \cite{BPS}, we rather require an incomplete condition satisfied by the elements of $\alpha(\HJ)$. Formally, the complex $\F_2\ra E\ra J[2]$ leads to the complex $\HJ \overset{\alpha}{\ra} \Hmu \overset{\mathbf{N}}{\ra} \Q^*/\Q^{*,2}$. Therefore we have the inclusion $\im(\alpha) \subset \ker(\mathbf{N})$; with no guaranteed equality. But this time, the morphism $\mathbf{N}$ is computable when restricted to $\etale \subset \Hmu$ as it is the canonical norm map. So we remove the elements of $\sLS \supset \SelT$ which are not in $\ker(\mathbf{N})$.

\subsubsection{Computing the set of bad primes}

In order to compute $\sLS$, we first need to figure out the primes of bad reduction for the Jacobian. It also guide the computation of $\SelF$ as primes of bad reduction are often better candidates for $\Tset$ than the ones of good reduction.\par
In general we compute a discriminant of the curve by double application of the resultant on the equation of the curve and rule out one by one the factors that are of good reduction. This gives a set that contains $\Sset$ since primes of bad reduction for the Jacobian are of bad reduction for the curve. We sometimes apply changes of variables when the computed discriminant is too huge to factor.\par
When $g=3$, one can do as in \cite{BPS} (See section 12.5).\par
Actually, according to \cite[Lem 7.1]{BPS}, we can replace $\sLS$ by $\widetilde{L(2,\mathcal{U})}$ for any set $\mathcal{U}$ that contains $2$, $\infty$ and the primes $p$ such that the Tamagawa number $c_p(J)$ of $J$ at $p$ is coprime to $2$. Unfortunately, computing $c_p(J)$ is quite a hard task in general.

\subsubsection{Computing Galois action and cohomology kernels explicitly}\label{section explicit Galois}

In a work to appear \cite{SHA}, we study the global and local cohomology kernels $\kF(k)$ in a theoretical, computational and statistical fashion. In particular, we describe how to explicitly compute the Galois action over the TC and over the various spaces mentioned in the present article from the image of the Galois representation $\rpztors : G_k \ra \Sp$. Computations and Galois action in $\E$ and $\Rw$ are handled. See also \cite[section 12.3.]{BPS}.\par
We use this to explicitly compute the space $\Sha^1(\Q,J[2])$ for a given curve. It is actually related to the global cohomology kernel $\kF(\Q)$ which only depends on the Galois group $G$. We therefore establish the data table of dimension of $\kF(G)$ for every conjugacy class of group $G\subset\Sp$ for $g=2,3,4$. We draw qualitative consequences from the observation of these data tables.\par
We show that the use of both Galois representations $\rpztors$ and $\rpzTC$ rather than one is more efficient to determine the global and local Galois groups of $J[2]$. Evidences for this are given in the Examples section of the present article.

\section{TC-descent}

We propose in this section to perform 2-descent with $\Delta := \Delta_{TC} $ being the whole set of Theta Characteristics rather than taking $\Delta:=\OTC$ as BPS do in \cite{BPS}. We call it \textit{TC-descent}. This avatar of 2-descent, entirely due to the author, have been implemented and used in practice. Its efficiency is discussed in the Examples section below.\par
Compared to OTC-descent, TC-descent always reduces the kernels and even sometimes releases us from the need for knowledge of the Galois group $G$ of the 2-torsion. This is all the more convenient since the latter is tricky to obtain in practice and stifles the automatization of the 2-descent. 
In a work to appear, we present theoretical and practicable conditions for the cohomology kernels to be trivial. 
The method could be useful for hyperelliptic curves as well, as long as the ETC Galois orbits are small enough. In particular, this could supplement the even genus case famously considered by Stoll in \cite{StollCH}.\bs

We explain how we construct the TC-descent map and discuss the advantages of having an injective descent map.


 

\subsection{Fake Descent Setup}\label{section Fake Descent setup}

Under the assumption of a rational point that we require, every class in $\Pic_X(K)$ can be represented by a divisor of $\DivOp_X(K)$. We can therefore find a rational divisor $\Sigma \in \DivOp(X\times \TC)(\Q)$ such that each $\Sigma_\beta$ with $\beta\in\TC$ is a divisor representing the corresponding TC and defined over its definition field. We thus obtain a \textit{fake descent setup} $(2,\TC,\Sigma)$ in the meaning of \cite[Def 6.7]{BPS} (with $D=K_X$ a fixed canonical divisor). In particular, the descent map $C_{TC} : J(k) \ra \etalek$ whose components are the Theta Fractions $(f_\beta)_{\beta\in \TC}$ coincides with the Cohomology map $\alpha_k : \HJk \ra \etalek$ restrained over $J(k)$ (See \cite[Prop 6.10]{BPS}).\par
Since $\TC $ is isomorphic to $\OTC \sqcup \ETC$ as a Galois-set, we observe $L_{TC}\simeq L_{Odd} \times L_{Even}$ where $L_{Even}$ is the \'etale algebra attached to $\ETC$. Loosely speaking, the proposed descent map is the BPS classic one extended with the even Theta Fractions.\par
We call this descent a \textit{TC-descent} in the present article.\bs

More generally, we can take $\Delta$ to be the union of any chosen orbits of TC, as long as the corresponding morphism $\pi_\Delta: E:=(\FDelta)_0 \ra J[2]$ is surjective (See section \ref{section Descent Cohomology}). Under the assumption of the existence of a rational point on the curve, this also defines a Fake Descent Setup. We call this descent an \textit{incomplete TC-descent} in the present article. It sometimes offers a sweet trade-off between efficiency and computation.

\subsection{The TC-descent map $C_{\TC}$ in practice}


As explained in \cite{TitiOTC}, due to the internal structure of the algorithm, it is not much more time consuming to return TC-descent data over OTC-descent data.\par
Moreover, the sets $\OTC$ and $\ETC$ are respectively of comparable sizes $2^{g-1}(2^g-1)$ and $2^{g-1}(2^g+1)$. As the Galois action over $\ETC$ is completely determined by the one over $\OTC$, it is reasonable to think that often when the orbits of one of the two sets split so do the orbits of the other. Obtaining the class groups and units groups of the fields of the \'etale algebra $L$ is the computational limiting factor of the 2-descent. We are therefore confident in the fact that the possibility of performing 2-descent with $\TC$ occurs quite often when it is with $\OTC$.

\subsection{Reducing the kernels}\label{section reducing kernel}

The local and global descent maps $C_{TC,k} : J(k)/2J(k) \ra L_{TC}^*/L_{TC}^{*,2}k^*$ associated with $\TC$ factor the usual descent maps $C_{OTC,k}$ used in \cite{BPS} derived from the Galois-set $\OTC$. Indeed, the commutative diagram 
\begin{tikzcd}
    E_{\OTC} \arrow[r] \arrow[rr, bend left] & E_{\TC}:= E_{\OTC} \oplus E_{\ETC} \arrow[r] & J[2]
\end{tikzcd}
leads to the commutative diagram
\begin{tikzcd}
    \Hop^1(k,J[2]) \arrow[r, "\alpha_{k,TC}"] \arrow[rr,bend left, "\alpha_{k,Odd}"] & \Hop^1(k,\mu_2^{\TC}/\mu_2) \arrow[r] & \Hop^1(k,\mT)
\end{tikzcd}
for any $k=\Q,\Q_p$. The kernels of the maps $C_{TC,k}$ are therefore always smaller (or equal) than the ones of $C_{OTC,k}$. It is not that surprising as in practice we just added the components $(f_\beta)_{\beta \in \ETC}$ attached to the even Theta Characteristics. As we soon show, there is a computational gap between trivial descent kernels and even small but non-trivial ones. Reducing the kernels is therefore of stark interest as trivial descent kernels greatly ease 2-descent in practice and that they are already observed rather small with the usual maps $C_{OTC,k}$.

\subsubsection*{Locally}
A non-trivial local descent kernel $\ker(C_p)$ forces us to carry out "the residue discs method" (See \cite{BPS} section 11.1.) which may be long to achieve. For each prime $p$ in the chosen set $\mathcal{T}$, the regular process consists in computing several different points of the local image $C_p(J(\Q_p)/2J(\Q_p))$ built from points in $X(\Q_p)$ in hopes of reaching the expected known dimension of $J(\Q_p)/2J(\Q_p)$. This usually works well in practice but is conditional on the injectivity of the map $C_p$ which is by the way ensured if successful. If the right dimension is not reached but we know that $C_p$ is injective (for instance because we have computed the local cohomology kernel $\kF(\Q_p)$ with help of the knowledge of the Galois group $G$ of the 2-torsion), then we know that we have to look for points of $J(\Q_p)$ expressed by divisors supported by points of the curves defined over some extensions of $\Q_p$. It occurs quite rarely. But if $C_p$ is not even injective, we are forced to do the same for one point in \textit{every} residue disc of small enough radius of $X(K)$ for \textit{every} extension $K$ of $\Q_p$ of a given known degree (a finite number by $p$-adic Kummer theory). See \cite[Lem 11.2.]{BPS}. In practise, the more the genus increases, the more complicated it becomes.\par
We hope for trivial local descent kernels but notice that, in opposite to the global situation, it provides a small computational asset to certify it in advance. In particular, computing $\dim_{\F_2} \kF(\Q_p)$ is mainly useful to construct $W_p$ in order to lower the global bound on the rank from $\SelF\times \kF(\Q)$ to $\SelF \times \sKF$. Moreover, if it does not fail, undertaking the regular process for computing $\im(C_p)$ proves local injectivity of $C_p$.

\subsubsection*{Globally}
Moreover, once we have computed all the local images, we then need to bound the dimension of the global descent kernel $\sKF$. 
It requires an explicit description of the Galois group $G\subset\Sp$ of $J[2]$ in order to compute $\kF(\Q)$. We can further reduce with use of the explicit knowledge of the local cohomology kernels $\kF(\Q_p)$ (which can turn out not to be trivial even if the computation of the local images showed that $C_p$ is injective). Therefore, if we can ensure that $\alpha_\Q$ is injective, it removes the need of $G$ in addition with sparing us the potential computation of the $\kF(\Q_p)$.\par
We study both theoretical and computational conditions on the Galois group for the TC-descent map to be injective in \cite{SHA}. We also list in \cite{SHA} the structure and TC/$J[2]$-orbits decomposition of all unfortunate groups $G\subset\Sp$ that may lead to a non-injective map $\alpha_\Q$ for $g=2,3,4$. In particular, partial knowledge (cardinal, group structure) of $G$ may be sufficient to skip the computation of the whole group description. Some of this ruling out partial knowledge may be acquired by local computations. Note as well that the orbits decomposition may be enough to guarantee the global (or local) injectivity of $C_{TC}$. The way we obtain it does not even require information about $G$ at all. See \cite{SHA} for more details.

\subsubsection*{Sum-up}
\begin{itemize}
    \item \textit{(Locally)} If  $\alpha_p$ is injective:\\
    Then $C_p$ is injective. Loosely speaking it means that we will probably be able to compute the local image $C_p(J(\Q_p)/2J(\Q_p))$ quite easily. Moreover we know that $W_p=0$ and we therefore do not need to know the local Galois group $D_p$.
    \item \textit{(Globally)} If $\alpha_\Q$ is injective:\\
    The global descent kernel $\sKF=0$ is trivial. There is therefore no need to know the Galois group $G$ and decompositions groups $D_p$ nor the spaces $W_p$.
\end{itemize}

\begin{rem}
    When the descent map $C_k$ is injective, we can also use it to detect linear independence in $J(k)/2J(k)$. This can be useful to exhibit finite index subgroups of $J(\Q_p)$ in view of applications \textit{\`a la} Chabauty-Coleman.
\end{rem}

\section{Lower bounds on the rank}

We present here few techniques to produce lower bounds of the rank. This sometimes allows us to certify the rank when it turns to match the upper bound obtained from descent.\bs

We denote $J_{c}(\Q):=\<[P-P_0] ~ |~P\in X(\Q),~ P\text{ known}\rangle$ the Jacobian subgroup spanned by differences of the rational points of the curve we know. 
We can handle the points $P$ that are singular for the chosen model by looking at the branches above $P$. For poles of the TF, we can move the canonical divisor $K_0$ by division by a well-chosen (or actually random) function of $\Hop^1(K_0)(\Q)$ as this does not change the image of $P$ under $C_k$.

\subsection{By the positive rank trick}\label{sectionPositiveRankTrick}

We present a method due to the author that allows us to sometimes certify that $\rank_\Q(J) =1$ when we have already proved that $\Sel(J) \leq 1$ (See \cite{PositiveRankTrick}).\par

In short, in the case like the one we assume where the curve possesses a rational point, at least one of the three following statements is true:
\begin{enumerate}
    \item $\rank_\Q(J) \geq 1$.
    \item $J[2](\Q)$ non-trivial.
    \item $\TC(\Q) \neq \varnothing$.
\end{enumerate}

The algorithm of the author \cite{TitiOTC} offers the possibility of computing the \'etale algebras $L_{J[2]}$, $L_{TC}$ and therefore of detecting potential rational TC and rational 2-torsion. 
Whether the statements 2. and 3. are true or not can therefore be checked quickly with what we need to compute anyway. This way, statement 1. is sometimes guaranteed for free. We give an example in the Examples section (See example \ref{example PRT}).

\subsection{By descent maps}\label{section lower bound map}

The maps we are about to present find a purpose in bounding the rank by below but can also be used to certify linear independence in $J(k)/2J(k)$ and a fortiori in $J(k)$ for both a global or a local field $k$.

\subsubsection{Global maps}

When the global descent map $C_\Q$ is injective, the dimension $\dim_{\F_2} C_\Q(J_c(\Q))$ gives a lower bound of the rank that may be sharp. We saw that TC-descent is always more likely to produce injective maps.\bs

\textit{Computing the global descent map $C_\Q$}\\
The main difficulty is that the global descent map is more complicated to compute than the local ones. Indeed, it breaks down to compute $C_\Q$ as the composition of two maps $C_\Q:J(\Q) \overset{\ev}{\ra} L^* \overset{\pi}{\ra }\etale$. The first one $\ev$ simply consists in the specialisation of the TF to a given divisor $D$ representing a rational class in the Jacobian and raises no issues in practice. This returns an element $\alpha:=\ev([D])\in L^*$ that we have to reduce into $\etale$. But, in opposite to the local case, no computable canonical representative stands out and squares are less easy to determine. Surely descent already led us to compute the space $\LS$ that contains a representative of the projection $\pi(\alpha)$. Moreover, the map $\pi$ is effective over $\LS$. On the other hand, it may be tricky to detect and remove the squares that relate $\alpha$ to some representative in $\LS$.\bs

\textit{Getting help from the norms: the Norm descent map $\eta$}\\
For this reason, it may sometimes be worth looking at the algebraic norm of the components of $\alpha\in L^*$. Let us consider the morphism of groups $\Nc: \alpha\in L^* \ra (\operatorname{Norm}_\Q(\alpha_i))_{i \leq w} \in (\Q^*)^w$ where $w$ is the number of fields constituting $L$. Further reducing, we obtain the morphism $\pi\circ \Nc : L^* \ra \NQ{w}$ where the quotient of $(\Q^*/\Q^{*,2})^w$ is by the diagonal action of $\Q^*$. We try to factor $\pi\circ \Nc$ through $\etale$. Note that $\Nc(L^{*,2}) \subset \Q^{*,2}$ and therefore $L^{*,2}$ is in the kernel of $\pi\circ \Nc$. But for $\lambda\in \Q^*$, we have $\pi\circ\Nc(\lambda) = (\lambda^{\epsilon_i})_{i\leq w}$ where $\epsilon_i\equiv \deg(L_i) [2]$ is the parity of the degree of the corresponding field of the \'etale algebra. Therefore, we rather split the fields in $L$ of even and odd-degree and consider the product of morphisms $\Nc: L^* \ra \Gamma:=\NQo \times \NQe$, where $w_o,w_e$ are respectively the number of odd and even-degree fields in $L$. Notice the asymmetry of the codomain $\Gamma$. It factors through $\etale$ indeed. We call the resulting factored morphism $\Nc_{etale}$. The commutative diagram below describes the situation.
\begin{equation*}
\begin{tikzcd}
    L^* \arrow[r, "\pi"] \arrow[d,"\Nc"] & \etale \arrow[d, dashrightarrow, "\Nc_{etale}"]\\
    (\Q^*)^{w_o} \times (\Q^*)^{w_e} \arrow[r, "\pi"] & \Gamma:= \NQo \times \NQe
\end{tikzcd}
\end{equation*}

Therefore, we are now equipped with a well-defined morphism $\eta:=\Nc_{etale}\circ C_\Q :J(\Q) \ra \Gamma$ that factor through $J(\Q)/2J(\Q)$. We call it the \textit{Norm descent map} in the present article. We also give the cohomological interpretation of $\Nc_{etale}$.

\begin{prop}
    Let $k$ be a number field or a $p$-adic field, $w_e$ and $w_o$ be the number of respectively even and odd-degree fields in the \'etale algebra $L$.\\
    There exists a well-defined morphism $\Nc_{etale} : \etalek \ra \Gamma_k:=(k^*/k^{*,2})^{w_o}/k^* \times (k^*/k^{*,2})^{w_e}$ given by taking the corresponding algebraic norms of the fields constituting $L$.\\
    Moreover, the morphism $\Nc_{etale}$ can be obtained by dualising and taking the first cohomology of the morphism of Galois-modules $(\F_2^{w_o})_0\times\F_2^{w_e} \ra E$.
\end{prop}

\begin{proof}
    We refer to the discussion above for the natural definition of $\Nc_{etale}$.\\
    For every TC orbit $\nu \subset \Delta$, we set $\mathbf{1}_\nu:=\sum_{\beta \in \nu} e_\beta \in \F_2^\Delta$. It comes with a morphism $\F_2 \ra \F_2^\Delta$ that maps $1$ to $\mathbf{1}_\nu$. If $\nu$ has even cardinal, $\mathbf{1}_\nu \in E$, and following the approach of \cite[Lem 6.16]{BPS}, by dualising and taking cohomology of $\F_2 \ra E$, we obtain a morphism $\Hop^1(k,\E)\simeq \etalek \ra \Hop^1(k,\mu_2)\simeq k^*/k^{*,2}$. We now take the product of the morphisms from all the orbits of even cardinal and consider the Galois-module $\F_2^{w_e}$ with trivial action. Dualising and taking the first cohomology of the morphism $e_\nu \in \F_2^{w_e} \mapsto \mathbf{1}_\nu \in E$, we obtain the second component $\Hop^1(k,\E)\simeq \etalek \ra \Hop^1(k,\mu_2^{w_e})\simeq \Hop^1(k,\mu_2)^{w_e} \simeq (k^*/k^{*,2})^{w_e}$ of the morphism $\Nc_{etale}$.\\
    The image of the equivalent map $\F_2^{w_o} \ra \F_2^\Delta$ is not contained in $E$. But restricting to the elements with an even number of components $(\F_2^{w_o})_0 \ra E$ gives the desired first component, $\Hop^1(k,\E) \ra \Hop^1(k,\mu_2^{w_o}/\mu_2)\simeq (k^*/k^{*,2})^{w_o}/k^* $ of $\Nc_{etale}$ once we dualize and take the first cohomology. We recall that $(k^*)^{w_o}$ is the \'etale algebra $l$ of the trivial Galois-set $\Delta:=[1,...,w_o]$ and that $\Hop^1(k,\mu_2^{w_o}/\mu_2) \simeq l^*/l^{*,2}k^*\simeq (k^*/k^{*,2})^{w_o}/k^*$ is the quotient of $(k^*/k^{*,2})^{w_o}$ by the diagonal action of $k^*$.
\end{proof}

For $k=\Q$, the Norm descent map $\eta$ can be easily computed in practice as the morphism $\pi \circ \Nc\circ \ev$. Indeed, reducing $\Nc(\alpha)= \Nc\circ\ev (D)\in (\Q^*)^{w_o} \times (\Q^*)^{w_e}$ into $\Gamma$ is now convenient as we are able to factor in $\Q$. Moreover, given a finite generating set of $J_c(\Q)$, we can perform linear algebra in $\Gamma$ gathering the primes appearing in the factorisation of their image by $\Nc\circ\ev$. Actually, up to squares and multiplications by elements of $\Q^*$, every tuple $\Nc\circ \ev(D)$ for $D\in J(\Q)$ can be represented in $\Gamma$ by the primes in $\Sset$ since $\SelF\subset (\etale)_{\Sset}$ (See \cite[Prop 7.2 + Theo 10.9]{BPS}). To ease computation, we can find some help in the fact that every element in $\Gamma$ has a canonical representative in $(\Q^*)^{w_o} \times (\Q^*)^{w_e}$ that has squarefree components and its smallest odd-degree field component equal to $1$.\bs

This allows us to ensure linear independence and non-triviality in $J(\Q)/2J(\Q)$ (and in $J(\Q)$) more easily that we would for $C_\Q$, at the expense however that $\eta$ may have a bigger kernel.

\begin{prop}
    If $J[2](\Q)=0$, then $\eta$ is trivial.
\end{prop}

\begin{proof}
    The image of the morphism $n: e_\nu \in(\F_2^{w_o})_0\times\F_2^{w_e} \mapsto \mathbf{1}_\nu \in E$ (see the proof above) is contained in $E(\Q)$ as the Galois-module $(\F_2^{w_o})_0\times\F_2^{w_e}$ has trivial action. Therefore the image of the composition $\pi\circ n : (\F_2^{w_o})_0\times\F_2^{w_e} \ra J[2]$ is contained in $J[2](\Q) = 0$. In other words, we have the inclusion $\im(n) \subset R$. According to the last property, dualising and taking cohomology of the complex $ (\F_2^{w_o})_0\times\F_2^{w_e} \overset{n}{\ra} E \overset{\pi}{\ra} J[2]$ leads to the complex $\Hop^1(\Q,J[2]) \overset{\alpha_\Q}{\ra} \Hop^1(\Q,\E) \overset{\Nc_{etale}}{\ra} \Gamma$. In particular, $\eta=\Nc_{etale}\circ C_\Q$ is trivial.
\end{proof}

Note however that in the unfortunate case where $\eta$ is trivial, we can use $\Nc$ instead of $\mathrm{N}$ to further reduce $\SelF$; See section \ref{section ST}.\bs

\textit{Non-trivial kernel}\\
We may still cope with a small non-trivial global descent kernel $\sKF$. Assume we have found a point $P\in J(\Q)$ with $C_\Q(P)=0$ (or $\Nc_{etale}\circ C_\Q(P)=0$). We want to prove that it is non-zero in $J(\Q)/2J(\Q)$. This would increase the lower bound given by $\dim_{\F_2} C_\Q(J_c(\Q))$.\par
It may be enough to check locally the vanishing of $C_p(P)$ at several primes $p$. If the local descent maps are injective and the Jacobian $J$ satisfies the local-global principle for the $2$-divisibility problem, we know the process will eventually terminate when $P$ is not a rational double. The author addresses the $2$-divisibility problem for Jacobians for its own sake in \cite{SHA}. Otherwise, slightly adapting the algorithm of \cite{Mascot} by shifting $0$ to $P$, one could compute a polynomial $\chi\in \Q[y]$ realising the \'etale algebra of $[2]^{-1}(P)$. The point $P\in J(\Q)$ is a rational double if and only if the polynomial $\chi$ has a linear factor.

\subsubsection{Local maps}
We can also perform a local equivalent by checking the dimensions $d_p:=\dim_{\F_2} C_p(J_c(\Q))$ for several primes $p$. If the restriction $J(\Q)/2J(\Q) \ra J(\Q_p)/2J(\Q_p)$ is injective for some $p$, we may be able to obtain a sharp lower bound. We are always guaranteed that $\rank_\Q(J) + \dim_{\F_2} J[2](\Q) \geq d_p$.\par
Moreover, if the restriction at $p$ is injective over $\sLS$, which is not that rare with $p=2$ for instance, and we verify the equality $d_p = \dim_{\F_2}(J[2](\Q_p))$, we are then guaranteed that $\rank_\Q(J) + \dim_{\F_2} J[2](\Q) = d_p$. In that case, we have certified the rank.\bs

Generalizing the approach of \cite[section 2.2]{IndependancyJQ} to our settings, we can combine the information gathered at multiple primes. Let $\mathcal{U}$ be a set of primes. We consider the map $\prod_{p\in \mathcal{U}} C_p: J_c(\Q) \ra \prod_{p\in \mathcal{U}} \Q_p^*/\Q_p^{*,2}$. Similarly, we have $\rank_\Q(J) + \dim_{\F_2} J[2](\Q) \geq \im (\prod_{p\in \mathcal{U}} C_p)$. Using Makdisi's algorithms (See \cite{Makdisi1}, \cite{Makdisi2}), the discussion about linear independence over $J(\Q)$ in \cite[section 2.2]{IndependancyJQ} also generalize to Jacobians via the maps $C_p$.\par
We show successful examples of this idea in the Example section. 

\subsection{By computations in the Jacobian}

We can also compute a multiple $B$ of the size $T$ of the rational torsion of the Jacobian by computing the cardinal of the $\F_p$-point $J(\F_p)$ for (several) prime(s) $p$. Taking the greatest common divisor, one deduces a multiple $B$ of $T$. See \cite[section 5.]{Pooneng2} or \cite{VanBommelTorsion}. It remains to pick a rational point $P\in J_c(\Q)$ and to check whether $B\cdot P$ is trivial using Makdisi's algorithms (See \cite{Makdisi1}, \cite{Makdisi2}). By testing few points, we may be able to exhibit non-torsion points. Furthermore, by checking linear dependence relations, we may be able to certify a matching lower bound.\bs

Following \cite[Ex 12.9.3.]{BPS}, we can choose a finite set of primes of good reduction $\Uset$ and look at the image of the morphism $J_c(\Q) \ra \prod_{p\in \Uset} J(\F_p)$. It may be easier to consider the morphism $\prod_{p\in \Uset} C_p$. If the $(C_p)_{p\in \Uset}$ are injective, it is as efficient. Moreover, primes of bad reduction are allowed.






\section{Examples}

In order to build trust and assess correctness of the algorithm, we present several 2-descents that we ran on already known examples of non-hyperelliptic curves and discuss them. The three first examples are the three first examples of the genus 3 non-hyperelliptic curves that Bruin, Poonen and Stoll considered on their original article \cite{BPS}. 
We could not apply our methods to the fourth curve as it has no rational points. We rederive the results about the curve \textit{$XE_2$} in \cite{Lombardo} and about a genus 4 non-hyperelliptic curve that Fu and Stoll computed in \cite{StollFu} with TC-descent. We perform an incomplete (but successful) TC-descent on a genus $5$ non-hyperelliptic curve in the next chapter.\bs

For each of the showcased examples, unless explicitly mentioned, nothing but an equation of the curve, a rational point and (sometimes) the places of bad reduction and/or the knowledge of the Galois group has been used to obtain the results. When it was necessary, we computed the Galois groups from their index in $\Sp$, from their cardinal or from an equation of the corresponding \'etale algebra using the online free Magma \cite{Magma} calculator. The time bottleneck of these examples lies in the restriction $\etale \ra L_p^*/L_p^{*,2}\Q_p^*$ at $p=2$. Indeed, the elements of the basis of $\Sset$-units spanning $\sLS$ can be quite involved and the restriction at $p=2$ is the most time consuming.
All computations have been run on the personal computer of the author.
The code has been written in the gp framework of the PARI/GP language \cite{PARI2} and is not optimized.\bs

In opposite to the previous methods, the algorithm of the author needs no Groebner basis computation. It is also automatized and may work over other (maybe higher but moderate) genus curves under the condition that the (Odd) TC orbits split enough so that we can compute the classgroups, the units groups of their corresponding number fields along with its global Galois group. The author would be happy to try performing 2-descent on other non-hyperelliptic curves. The reader should not hesitate to contact him with a given curve (of moderate genus) to perform descent along with a short motivation of the original purpose.

\subsection{A plane quartic}

We perform TC-descent on the first example of \cite[section 12.9.1.]{BPS}. Notice that the global cohomology and descent maps are injective while the OTC ones are not. As for genus 3 curves, rather than taking $\Sset$ to be the set of places of bad reduction, we can refine and take only places of trivial Tamagawa number by looking at the discriminant (See \cite[section 12.6.3.]{BPS}). In this example, it reduces $\Sset$ from $\{2,29,163 \}$ to $\{ 2 \}$.\bs

At $p=2$, one can only rule out all but 3 possible local Galois groups by looking at the permutation type of $\rpzOTC(\Frob_p)$. Looking at the Frobenius type of $\rpzTC(\Frob_p)$ and $\rpztors(\Frob_p)$, we managed to distinguish the right group. See the section \ref{section explicit Galois}.

\subsubsection*{Technical outputs}
wall time computation (once (O)TC and Galois groups computed): 3min, 13,637 ms\\
Curve: $yx^3 + (-y^2 - 1)x^2 + (-y^2 + 1)x + y^3=0$,		$P_0 = [0, 1, 0]$, $g = 3$\\
$\Sset = \{2\}$, $\Tset = \{2\}$, $\Delta = \TC$\\
Rational Theta Fraction: $-yx + (y^2 + z^2)$\\
Orbit Decomposition. $\OTC: 28_1$   $\ETC: 1_1 35_1$\\ 
Orbit Decomposition. $J[2]\setminus \{0\}: 28_1 35_1$\\ 
$\dim_{\F_2} J[2](\Q) = 0$\\
$\dim$ classgroups: 1, 1, 1\\   
$\dim_{\F_2} \ker(\alpha_Q)= 0$,   Global cohomology map $\alpha_\Q$ injective? yes\\
$\dim H^1_{\Sset}(\Q,\E)= 38$,    $\dim_{\F_2} H^1_{\Sset}(\Q,\E)_{\text{norm}}= 38$\\

$p = 2$\\
        Orbit Decomposition. $\OTC: 28_1$   $\ETC: 1_1 7_1 28_1$\\
        $\dim_{\F_2} \JQp= 3$,   $\dim_{\F_2} \im(C_p)= 3$,   complete? yes\\
        $\dim_{\F_2} \ker(\alpha_p)= 0$,   $\dim_{\F_2} \ker(C_p)= 0$,   local descent map $C_p$ injective? yes\\
        $\dim_{\F_2} W_p= 0$\\
        $\dim_{\F_2} S_{\{2\} } = 0$\\
        
$\coker(\kappa)$ trivial? yes\\
$\Sha^1(\Q,\Rw)$ trivial? yes\\
$\dim_{\F_2} \SelF\leq 0$,   $\dim_{\F_2} \sKF\leq 0$,   $\dim_{\F_2} \Sel\leq 0$,  $\rank_\Q(J)\leq 0$\\
result from the original article: $\rank_\Q(J) = 0$

\subsection{An injective TC-descent of certified rank 1}\label{example PRT}

We perform TC-descent on the second example of \cite[section 12.9.2.]{BPS}. Notice that the global cohomology map is injective. We can certify the rank using the positive rank trick (See section \ref{sectionPositiveRankTrick}) and/or the local descent maps (See section \ref{section lower bound map}).

\subsubsection*{Technical outputs}

wall time computation (once (O)TC and Galois groups computed): 15min, 6,113 ms\\
Curve: $-yx^3 + (y^2 + 1)x^2 + x + (-y^3 - 2y^2 - y)=0$,		$P_0 = [1, 0, 0]$, $g = 3$\\
Known rational points: $\{[0, 1, 0], [1, 0, 0], [1, 1, 0], [0, -1, 1], [0, 0, 1]\}$\\
$\Sset =\{2\}$, $\Tset = \{2 \}$, $\Delta = \TC$\\
No rational TC found\\
Orbit Decomposition. $\OTC: 28_1$  $\ETC: 36_1$\\ 
Orbit Decomposition. $J[2]\setminus \{0\}: 63_1$\\ 
$\dim_{\F_2} J[2](\Q) = 0$\\
$\dim$ classgroups: 1, 1\\   
$\dim_{\F_2} ker(\alpha_\Q)=0$,    Global cohomology map $\alpha_\Q$ injective? yes\\
$\dim_{\F_2} H^1_{\Sset}(\Q,\E)= 37$,   $\dim_{\F_2} H^1_{\Sset}(\Q,\E)_{text{norm}}= 37$\\

$p = 2$\\
        Orbit Decomposition: $\OTC : 28_1$  $\ETC: 8_1 28_1$\\
        $\dim_{\F_2} \JQp =   3$,   $\dim_{\F_2} \im(C_p)= 3$,  complete? yes\\
        $\dim_{\F_2} \ker(\alpha_p)= 0$,   $\dim_{\F_2} \ker(C_p)= 0$,   local descent map $C_p$ injective? yes\\
        $\dim_{\F_2} W_p= 0$\\
        $\dim_{\F_2} S_{\{2 \}}=1$\\

        $\dim_{\F_2} C_p(J_c(\Q)) -\dim_{\F_2} J[2](\Q) = 1$\\
        $\leadsto \rank_\Q(J) \geq 1$\\

$\coker( \kappa)$ trivial? yes\\
$\dim_{\F_2} \Sha^1(\Q,\E)\leq 2$\\
$\dim_{\F_2} \SelF\leq 1$,	$\dim_{\F_2} \sKF\leq 0$,	$\dim_{\F_2} \Sel \leq 1$,	$\rank_\Q(J) \leq 1$\\
Positive Rank Trick: $ 1\leq \rank_\Q(J)$ \hspace{5mm} $\leadsto \rank_\Q(J)=1$\\
result from the original article: $\rank_\Q(J) = 1$\\

Notice that the positive rank trick (See section \ref{sectionPositiveRankTrick}) applies since $J[2](\Q)$ is trivial and the curve has no rational TC.

\subsection{TC-descent with primes of good reduction}

We perform TC-descent on the third example of \cite[section 12.9.3.]{BPS}. At $p=2$, the local cohomology kernel $\kF(\Q_2)$ has dimension $2$. Rather than using the "residue discs" method, we compute the local images at the primes of good reduction
$$p \in \Tset:=\{ 29, 31, 37, 167, 193 \}.$$
For the three first primes $p$, we observed that $\dim_{\F_2} J(\Q_p)/2J(\Q_p) = \dim_{\F_2} J[2](\Q_p) =0$. Weirdly enough, there is no rational $2$-torsion for almost all primes below $200$. To impose non-trivial conditions, we had to pick the last two primes of $\Tset$. In both cases, the local descent map was injective and we managed to recover the full image. But we had to draw more local points than usual.\bs

However, $C_2$ is injective on $J_c(\Q)$ and the lower bound obtained that way matches the descent bound. We are therefore able to certify the exact rank.

\subsubsection*{Technical Outputs}

Curve: $yx^3 + (-2y^2 - y - 2)x^2 + (y^3 - y^2 + 2y + 3)x + (y^3 - y)=0$,                $P_0 = [1, 0, 0]$, $g = 3$\\
$\Sset = \{2, 13\}$, $\Tset = \{29, 31, 37, 167, 193\}$, $\Delta = \TC$\\
Orbit Decomposition. $\OTC: 28_1$  $\ETC: 36_1$\\ 
Orbit Decomposition. $J[2]: 63_1$\\ 
$\dim_{\F_2} J[2](\Q)= 0$\\
$\dim$ classgroups: 1,   1\\   
$\dim_{\F_2} \ker(\alpha_\Q)= 2$,    Global cohomology map $\alpha_\Q$ injective? No\\
$\dim_{\F_2} \Hop^1_\Sset(\Q,\E)= 45$,   $\dim_{\F_2} \Hop^1_\Sset(\Q,\E)_{\text{norm}}= 45$\\

$p = 29$\\
$\dim_{\F_2} J(\Q_p)/2J(\Q_p)= 0$,   $\dim_{\F_2} S_{\{29 \}}= 36$\par
$p = 31$\\
$\dim_{\F_2} J(\Q_p)/2J(\Q_p)=   0$,   $\dim_{\F_2} S_{\{29,31 \}}= 28$\par
$p = 37$\\
$\dim_{\F_2} J(\Q_p)/2J(\Q_p)=   0$,   $\dim_{\F_2} S_{\{29,31,,37 \}}= 21$\par
$p = 167$\\
$\dim_{\F_2} J(\Q_p)/2J(\Q_p)= 3$,   $\dim_{\F_2} \im(C_p) = 3$,   complete? yes\\
$\dim_{\F_2} S_{\{29,31,37,167 \}}= 2$\par
$p = 193$\\
$\dim_{\F_2} J(\Q_p)/2J(\Q_p)= 3$,   $\dim_{\F_2} \im(C_p) = 3$,   complete? yes\\
$\dim_{\F_2} S_{\{29,31,37,167,193 \}}= 1$\\

\noindent
$\dim_{\F_2} \SelF\leq 1$,	$\dim_{\F_2} \sKF\leq 2$,	$\dim_{\F_2} \Sel \leq 3$,	$\rank_\Q(J) \leq 3$\\
rational points (not complete): $[1, 1, 0], [0, 1, 0], [1, 0, 0], [0, -1, 1], [0, 0, 1], [0, 1, 1]$\\
$\dim_{\F_2} C_2(J_c(\Q)) = 3$ \hspace{5mm} $\leadsto \rank_\Q(J) =3$\\
result from the original article: $\rank_\Q(J) = 3$

\subsection{An injective TC-descent but non-injective OTC-descent}

We perform 2-descent on the curve $XE_2$ presented in \cite{Lombardo}. It is a genus 3 non-hyperelliptic curve. The authors of \cite{Lombardo} succeeds to perform 2-descent using special features of the curve, namely that it is a twist of the Klein Quartic. For such curves, one can often perform Fake descent with the "triangles" method \cite[Ex 6.9.iii.]{BPS} by choosing $\Delta$ to be a set of particular sums of Weierstrass points. We tried to perform OTC-descent on the curve but the local map turned out to be non-injective at $p=2$. We therefore could not certify the local image. However, we have been able to perform TC-descent that turned to be injective.\par
Notice that we have been compelled to compute the local images of several prime to obtain a sharp bound on the Selmer group since $\Tset=\{ 2 \}$ was not enough.\par
We also verified with success all the descents performed on the curves $(XE_i)_i$ in \cite{Lombardo} but do not display the technical outputs for space saving.

\subsubsection*{Technical outputs}

wall time computation (once (O)TC and Galois groups computed): 11min, 27,807 ms\\
Curve: $(-2y + 2)x^3 + 3x^2 + (-2y^3 - 6y + 2)x + (4y^4 - 12y^3 + 3y^2 + 7y - 6)=0$,		$P_0 = [1, 0, 0]$, $g = 3$\\
$\Sset = \{2, 3, 7\}$, $\Tset = \{2\}$, $\Delta = \TC$\\
Theta Fraction Rational TC: $-x^3 + (-y^2 - zy - 2z^2)x + (2y^3 - 4zy^2 - 2z^2y + 2z^3)$\\
Orbit Decomposition. $\OTC: 28_1$, $\ETC: 1_1 14_1 21_1$\\ 
Orbit Decomposition. $J[2]\setminus \{0\}: 14_1 21_1 28_1$\\
$\dim_{\F_2} J[2](\Q)= 0$\\
$\dim$ classgroups: 1, 1, 1, 3\\
$\dim_{\F_2} \ker(\alpha_\Q)=0$,  Global map $\alpha_\Q$ injective? yes\\
$\dim_{\F_2} H^1_{\Sset}(\Q,\E)= 62$,   $\dim_{\F_2} H^1_{\Sset}(\Q,\E)_{\text{norm}}= 62$\\

$p = 2$\\
    Orbit Decomposition. $\OTC: 4_1 12_2$,  $\ETC: 1_2 3_2 4_1 6_2 12_1$\\ 
    $\dim_{\F_2} \JQp = 4$,   $\dim_{\F_2} \im(C_p)= 4$,  complete? yes\\
    $\dim_{\F_2} \ker(\alpha_p)= 0$,   $\dim_{\F_2} \ker(C_p)= 0$,   local descent map $C_p$ injective? yes\\
    $\dim_{\F_2} W_p= 0$\\
    $\dim_{\F_2} S_{\{2\}}= 4$\\

    $p = 3$\\
Orbit Decomposition. $\OTC: 1_1 3_1 6_4$, $\ETC : 1_1 2_1 3_1 6_5$\\ 
$\dim_{\F_2} \JQp = 1$,    $\dim_{\F_2} \im(C_p)= 1$,   complete? yes\\
$\dim_{\F_2} \ker(\alpha_p)= 0$,    $\dim_{\F_2} \ker(C_p)= 0$,   local descent map $C_p$ injective? yes\\
$\dim_{\F_2} S_{\{2,3 \} } = 2$\\
    
$\coker(\kappa)$ trivial? yes\\
$\Sha^1(\Q,\Rw)$ trivial? yes\\
$\dim_{\F_2} \SelF \leq 2$,   $\dim_{\F_2} \sKF\leq 0$,	 $\dim_{\F_2} \Sel \leq 2$,	$\rank_\Q(J) \leq 2$\\
rational points (not complete): $[1, 1, 0], [1, 0, 0]$\\
$\dim_{\F_2} C_3(J_c(\Q)) =1$. So $1\leq \rank_\Q(J)$ \hspace{5mm} $\leadsto 1 \leq \rank_\Q(J) \leq 2$\\
result from the original article: $\rank_\Q(J) = 2$

\subsection{A failed attempt in genus 4}

We show an example of what could go wrong while trying to perform $2$-descent on a curve. In \cite{17T7}, a group of authors shows that \textit{17T7} is a Galois group over $\Q$. They relate it to the genus $g=4$ projective curve defined by the equation

\begin{minipage}{\dimexpr\textwidth-0cm}
\tiny
$8x^4y + 8x^3y^2 + 10x^2y^3 + 4xy^4 + 2y^5 - 8x^4 + 24x^3y + 6x^2y^2 + 12xy^3 - 2y^4 - 2x^3 - 3x^2y + 6xy^2 - 11y^3 + 14x^2 - 6xy - 4y^2 - 7x + 2y + 1=0.$
\end{minipage}

We tried to perform $2$-descent on. Unfortunately, our alorithm returns that the orbit decomposition of the TC is
$$ \OTC: 120_1 \hspace{5mm} \ETC: 136_1.$$
On this curve, $2$-descent is therefore out of computational reach as it would demand to compute the classgroup, units groups (and TF) of a number field of degree at least $120$.

\subsection{A genus 4 non-hyperelliptic curve}

We reproduce OTC-descent on the curve presented in \cite{StollFu}. Indeed, TC-descent involves a bigger etale algebra and therefore makes $\sLS$ bigger. It was not necessary as the OTC-descent turned to be already injective. In opposite to \cite{StollFu}, we obtain the Galois action on the OTC by the automatized algorithm \cite{TitiOTC}. As for the first example displayed in the present article, we managed to distinguish the right local group at $p=2$ by looking at the permutation type not only of $\rpzOTC(\Frob_p)$ but also of $\rpzTC(\Frob_p)$ and $\rpztors(\Frob_p)$.\par
We also exhibit a non-torsion point on $J(\Q)$, therefore certifying the rank exactly.

\subsubsection*{Technical outputs}

wall time computation (once (O)TC and Galois groups computed): 10min, 38,587 ms\\
Curve: $(y^2 - y)x^3 + (y^3 - y + 1)x^2 + (-y^3 - y^2 + 2y - 1)x + (y^2 - y)=0$,		$P_0 = [1, -1, 0]$,   $g = 4$\\
$\Sset = \{2\}$,   $\Tset = \{2 \}$,   $\Delta = \OTC$\\
Theta Fraction Rational TC: $yx - z^2$, $x+y$, $(2y - 1)x - y$\\
Orbit Decomposition. $\OTC: 1_3 3_3 12_1 24_4$, $\ETC: 1_1 3_1 12_3 24_4$\\ 
Orbit Decomposition. $J[2]\setminus \{0\}: 1_3 3_4 12_4 24_8$\\ 
$\dim_{\F_2} J[2](\Q)= 2$\\
classgroup: 1, 1, 1, 1, 1, 1, 1, 1, 1, 1, 1\\   
$\dim_{\F_2} \ker(\alpha_\Q)=0$,   Global map $\alpha_\Q$ injective? yes\\
$\dim_{\F_2} H^1_{\Sset}(\Q,\E)= 94$,   $\dim_{\F_2} H^1_{\Sset}(\Q,\E)_{\text{norm}}= 94$\\

$p = 2$\\
Orbit Decomposition. $\OTC: 1_8 2_6 6_6 8_8$\\
$\dim_{\F_2} \JQp = 8$,   $\dim_{\F_2} \im(C_p)= 8$,  complete? yes\\
$\dim_{\F_2} \ker(\alpha_p)= 0$,   $\dim_{\F_2} \ker(C_p)= 0$,   local descent map $C_p$ injective? yes\\
$\dim_{\F_2} W_p= 0$\\
$\dim_{\F_2} S_{\{2 \} } = 3$\\

$\coker(\kappa)$ trivial? yes\\
$\dim_{\F_2} \SelF \leq 3$, 	$\dim_{\F_2} \sKF\leq 0$,	$\dim_{\F_2} \Sel \leq 3$,    $\rank_\Q(J)\leq 1$\\
Checking torsion: $[0,0,1]-[1,-1,0]\in J(\Q)$ has infinite order \hspace{5mm} $\leadsto \rank_\Q(J) = 1$\\
result from the original article: $\rank_\Q(J) = 1$


\section{A more involved example:\\
A genus 5 non-hyperelliptic curve}

We performed TC-descent over the projective genus $5$ non-hyperelliptic curve $X$ given by the affine singular plane model
\begin{gather*}
\parbox[t]{1.0\displaywidth }{\raggedright{\small
$
    (6y^4 + 13y^3 + 18y^2 + 13y + 6)x^4 + (26y^5 + 79y^4 + 131y^3 + 131y^2 + 79y + 26)x^3 + (48y^6 + 190y^5 + 384y^4 + 484y^3 + 384y^2 + 190y + 48)x^2 + (45y^7 + 218y^6 + 520y^5 + 785y^4 + 785y^3 + 520y^2 + 218y + 45)x + (17y^8 + 94y^7 + 249y^6 + 422y^5 + 500y^4 + 422y^3 + 249y^2 + 94y + 17) =0.
$
}}
\end{gather*}

It is labelled \texttt{MC 67.68.5.a.1} in the LMFDB database. The curve is isomorphic to the modular curve $X_0(67)$. However, we did not make use of any of the modular features of the curve to perform $2$-descent on and ran the algorithm in its full generality. We picked a modular curve for three reasons. First, using the BSD conjecture, we could check our results by comparing with the analytic rank. Secondly, we could therefore choose a curve with a non-zero rank as we though this would produce a meaningful example. Finally, we expected the curve to be special enough for its orbits to split. In particular, we do not claim that we could perform $2$-descent on a general non-hyperelliptic genus $5$ curve. In fact the corresponding \'etale algebra and TF would be out of reasonable computational reach for our algorithm.\bs

To be more accurate, we performed an \textit{incomplete TC-descent} (See section \ref{section Fake Descent setup}). We computed the TC-orbits of the curve whose degree are $\leq 30$ using the author's algorithm \cite{TitiOTC}. This resulted in four OTC orbits of degree $1,15,30,30$ and a single degree $3$ ETC orbit. In fact, even if small, the ETC orbit turned to be essential for the morphism $\pi: E \ra J[2]$ to be surjective. From now on, we refer to these orbits as the \textit{small orbits}.\par
We were able to compute the \'etale algebra of all TC while selecting only small orbits for the computations of the TF.
This is possible with our method \cite{TitiOTC}, but not \textit{a priori}  with other methods based on Groebner bases, as those proceed by first computing all the Theta Fractions of the Odd TC only before projecting the equations over $\mathbb{P}^1$ to realise the \'etale algebra.
Moreover, the Groebner basis based methods cannot be \textit{a priori} adapted to ETC as they have no geometrical counterpart. In general, ETC can be expressed by the OTC as the latter span $J[2]$; however, in the case of this curve, the OTC small orbits generate neither the full space $J[2]$ nor the size $3$ ETC orbit. We recall the necessity of the ETC orbit for this specific curve. In particular, using a Groebner basis computation would force us to compute bigger OTC orbits than we have.
Regarding our method at least, the remaining orbits appeared to require more $p$-adic precision and therefore more time computation than what we used so far. We would expect the computation to be quite involved if one wanted to undertake it since the small orbits already demanded about a day and a half.\bs

Part of the computations were realised with Magma (See \cite{Magma}). The rest of the computations were realised in the PARI/GP software (See \cite{PARI2}).

\subsection{Computing the Galois groups}

The outputs of the algorithm \cite{TitiOTC} contains an explicit bijection $J[2] \simeq \F_2^{2g}$ along with an explicit dictionary $\{ TC\} \leftrightarrow \F_2^{2g}$ that is compatible under the sum. We were therefore able to check easily that the selected TC orbits span the whole group $J[2]$. This ensures that the $2$-descent is valid. But it also implies that the Galois representation attached to the TC factors through the one afforded by the small orbits. The \texttt{GaloisGroup} Magma function (very quickly) provided an equation for the Galois closure $F$ of the field associated to the unique orbit of TC of cardinal $45$. We checked (quite quickly as well) that the field of definition of the small orbits were included in $F$ using the \texttt{nfisincl} PARI/GP function (See the technical outputs). This identifies the Galois group $G$ of the TC as the group $\operatorname{Gal}(F/\Q)$.\par
In particular, we could determine that $G$ is of cardinal $360$, as well as its group structure (Magma label $S3*A5$). From the factorisation of the \'etale polynomial $\chi$ of the TC we also had the orbit decomposition. With Magma again, we were able to run through all the conjugacy classes of subgroups of $\Sp$ of cardinal $360$ and rule out all but 3 groups based on these features. It took few hours. Further looking at the orbit decomposition of the OTC and ETC separately we identified the right group $G$ as embedded in $\Sp$. For generators of a group in the $\Sp$-conjugacy class, see the technical outputs.\par
We also identified all the decomposition groups $D_p$ of $G$. Magma taught us that $G$ had only $11$ $G$-conjugacy classes of cyclic subgroups. Each of them was uniquely distinguished from the others by its TC orbit decomposition (with OTC and ETC separated). Therefore, by factoring the \'etale polynomial $\chi$ of the TC over $\F_p$ for many primes $p$ we recovered an explicit prime for each conjugacy class of Frobenius in minutes. We then certified the orbit decomposition by factoring $\chi$ over $\Q_p$. The same process turned to be successful as well for the two primes of bad reduction $2$ and $67$ of the curve.

\subsection{The 2-descent}

As we explained in the above action, we could identify an explicit prime $p$ for each conjugacy class of decomposition subgroup $D_p$ of $G$. This gives a set
$$ \Tset_{Cebo}:=\{\infty, 2,67,9439,277,23,293,1163,263,17,271,239,241 \}.$$
Note the two first primes $2,67$ are the unique bad primes of the curve. The primes are then ordered by increasing size of decomposition subgroup. For instance, $9439$ is a splitting prime. We set 
$$\Tset := \{\infty,67,9439,277,23,293,1163,263,17,271,239,241 \} $$
as the places for the local images. Note that $\Tset \cup \{2 \} = \Tset_{Cebo}$. Indeed, $2$ is the only place $p$ for which the local map $C_p$ is not injective. We reached a $3$-dimensional subspace out of the expected $5$-dimensional one and therefore cannot use this prime to impose descent condition. For all the other places, the local descent map $C_p$ is injective and we computed the full local image easily. This returns a bound for $\SelF$ of size
$$ \dim_{\F_2}S_{\Tset} = 3.$$
We call $I_2\subset\im(C_2)$ the subspace that we managed to compute. We have $\dim_{\F_2} I_2 =3$ and $\dim_{\F_2} \im(C_2) \leq \dim_{\F_2} J(\Q_2)/2J(\Q_2) =5$. We observed that the restriction $\res_2$ at $p=2$ is injective over $S_{\Tset}$ and that the intersection of the image $\res_2(S_{\Tset})$ with $I_2$ is trivial. Since we are not guaranteed that $I_2 =\im(C_2)$ based on the dimensions, we cannot certify that $S_{\Tset \cup \{2\}}=0$. However, with $U$ a supplement $I_2 \oplus U = \im(C_2)$, we can still ensure that $\dim_{\F_2} S_{\Tset \cup \{2\}} \leq \dim_{\F_2} U \leq \dim_{\F_2} \im(C_2)-3\leq 2$. Note that $S_{\Tset \cup \{2\}} = S_{\Tset_{Cebo}}= \SelF$ as the places in $\Tset_{Cebo}$ cover all the Frobenius (See \cite[Thm 10.9]{BPS}). So we have
$$ \dim_{\F_2} S_{\Tset_{Cebo}} = \dim_{\F_2} \SelF \leq 2.$$
From the explicit Galois group $G$, we managed to compute
$$ \dim_{\F_2} \kF(\Q) = 2.$$
Unfortunately, from the explicit decomposition subgroups $D_p$, we found out that all the local cohomology kernels $\kF(\Q_p) =0$ are trivial for $p\neq 2$. We therefore cannot reduce furthermore the bound $\kF(\Q)$ using the map $\kappa_{\Tset} : \kF(\Q) \ra \Pi_{\pi \in \Tset} W_p=0$ (See section \ref{sectionKernel}).\par
Factoring the \'etale polynomial $\chi$ of $J[2]$, we checked that $J[2](\Q)=0$ is trivial. Finally, we can certify the bounds
$$ \dim_{\F_2} \SelF \leq 2, \hspace{5mm} \dim_{\F2} \Sel \leq 2+2, \hspace{5mm} \rank_\Q(J) \leq 4.$$
According to the BSD conjecture, we have $\rank_\Q(J) = 2$; obtained from the analytic rank. Notice that we can ensure that the Chabauty-Coleman condition $\rank_\Q(J) < g$ is satisfied.\bs

We now present \textit{ad-hoc} techniques to improve the latter results.\par
We recall that $\sKF:=\ker(\alpha_{|\Sel})\simeq \ker(\kappa)$ where $\kappa : \kF(\Q) \ra \Pi_{p} W_p$ (See \cite[Thm 10.14.]{BPS} or section \ref{sectionKernel}). Therefore, in our case, $\sKF\simeq \ker (\kF(\Q) \ra \kF(\Q_2)/\ker(C_2))$. We checked that the restriction $\kF(\Q) \ra \kF(\Q_2)$ was a bijection. We give more details in how to perform this general computation in \cite{SHA}. In particular, $\dim_{\F_2} \sKF = \dim_{\F_2} \ker(C_2)$.\par
Now, we already mentioned that since the restriction at $p=2$ is injective over $S_{\Tset}$ and that the intersection of the image $\res_2(S_{\Tset})$ with $I_2$ is trivial, we have $\dim_{\F_2} \SelF \leq \dim_{\F_2} \im(C_2) -3$.\\
Gathering these two facts, we have that
$$ \rank_\Q(J) \leq \dim_{\F_2}\Sel \leq \dim_{\F_2} \sKF + \dim_{\F_2} \SelF \leq \dim_{\F_2} \ker(C_2) + \dim_{\F_2} \im(C_2) -3.$$
We recall that, from the rank theorem, $\dim_{\F_2} \ker(C_2) + \dim_{\F_2} \im(C_2) = \dim_{\F_2} J(\Q_2)/2J(\Q_2) = 5$. In other words, we have proven that
\begin{center}
    \fbox{$\rank_\Q(J) \leq \dim_{\F_2}\Sel \leq 2.$}
\end{center}

\subsection{Expectations}

We computed the local descent kernel $\kF(\Q_2)$ at $p=2$ and found out its dimension is $\dim_{\F_2} \kF(\Q_2)=2$. It is consistent with the fact that a $2$-dimensional space was missing from our computed attempt for the local image $\im(C_2)$. Or in other words that $\dim_{\F_2} J(\Q_2)/2J(\Q_2)-\dim_{\F_2} I_2=2$. Nevertheless, it is not certified that $\ker(C_2) = \kF(\Q_2)$ as we expect it to be. In order to prove that, one would have to perform the residue discs method (See section \ref{section reducing kernel}) at $p=2$. If one succeeds or assuming the local image at $p=2$ is of dimension $3$ as we found out, we would obtain $S_{\Tset_{Cebo}}=S_{\Tset \cup \{2\} } = 0$. This would lead to the expected bound $\rank_\Q(J) \leq \dim_{\F_2} \Sel \leq \dim_{\F_2} \SelF + \dim_{\F_2} \kF(\Q) \leq 0+2 $.\par
We computed the global cohomology kernel $\kF(\Q)$ and the local one $\kF(\Q_2)$ at $p=2$ for the full TC-descent $\Delta = \TC$. They are both of dimension $2$. Unless the corresponding local descent map $C_2$ is surprisingly injective over $J(\Q_2)/2J(\Q_2)$, it is unnecessary to compute all the TC to expect better results.

\subsection{Lower Bound}

We set $J_c(\Q):=\< [1, -3, 1]-[0,-1,1], [-1, 1, 0]-[0,-1,1] \ran \subset J(\Q)$. The trivial points $[-1, 0, 1]$ and $[0, 1, 0]$ are singular on the model with no rational branches.\par
As the $2$-torsion is trivial, we cannot use $\eta$ to prove non-triviality of rational points on the Jacobian.\par
We computed the image of $\Pi_{p\in \Tset_{Cebo}} C_p :J_c(\Q) \ra \Pi_{p\in \Tset_{Cebo}} L_p^*/L_p^{*,2}$ that turned out to be trivial. The only non-injective $C_p$ is at $p=2$. The cardinal $|J(\F_2)|$ is odd and thus every element in $J(\F_2)$ is a double. Therefore, to show non-triviality in $J(\Q)/2J(\Q)$ of a class $P\in J(\Q)$ by $p$-adic methods, one would have to find another model of the curve $X$ with good reduction at $p=2$ and approximate $P$ in $J(\Q_2)$.\par
Nevertheless, by computing a multiple of the cardinal of the rational torsion $J_{tors}(\Q)$, we were able to check that $D=[1, -3, 1]-[0,-1,1]\in J(\Q)$ is non-torsion. In particular, we get the bound
\begin{center}
    \fbox{$1\leq \rank_\Q(J).$}
\end{center}

Unfortunately, even if we found another explicit non-torsion point, we could not certify that they are linearly independent. Therefore, we could not improve the lower bound that way. Anyway, for the record, the point $[-1, 1, 0]-[0,-1,1]$ is $11$-torsion.

\subsection{A counter-example to the Local-Global Principle}

For any set $\Uset$ and any Galois-module $M$ over $k$, we denote $\Sha^1_{\Uset}(k,M)$ the Shafarevich-like group of $M$ with trivial conditions imposed everywhere but in the places of $\Uset$. Since $\kF(\Q_p)=0$ for all $p\in \Tset$ and that this latter set covers all the possible Frobenius and all but $p=2$ of the bad places, we have $\kF(\Q) \subset \Sha^1_{ \{2\} }(\Q,J[2])$. In particular,
$$ \dim_{\F_2} \Sha^1_{\{ 2\}}(\Q,J[2]) \geq 2.$$
Moreover, the images of $J_c(\Q) $ under $C_p$ for each $p\in \Tset$ are trivial and the corresponding maps $\alpha_p$ are injective, so 
$$J_c(\Q) \subset \Sha^1_{ \{2\} }(\Q,J[2])$$ 
as well. We showed that the point $D=[1, -3, 1]-[0,-1,1]\in J(\Q)$ is of infinite order but we could not guarantee that $D$ is not a double. Doing so would prove that $D$ is a counterexample to the Local-Global principle for the $2$-divisibility problem in $J$.\bs

In a work to appear, we present methods to systematically compute $\Sha^1(\Q,J[2])$ from the Galois groups $G$ of $J[2]$ and its decomposition subgroups at primes of bad reduction. It does not involve computations in $\etale$. In particular, we study more carefully the Shafarevich-like groups $\Sha^1(\Q,J[2])$ and $\Sha^1_{\{2 \} }(\Q,J[2])$ of $X$.\par


\subsection*{Technical outputs}

The author has gathered data related to this TC-descent on his personal public GitHub repository
\texttt{ThibautMisme/ArticleDescentExamples}.\\
It includes a description of the Galois group $G$, short equations for the fields of the small orbits, the \'etale polynomials $\chi$ for the \'etale algebras and the TF of the small orbits.\\
Please visit \url{https://github.com/ThibautMisme/ArticleDescentExamples}.\bs

\noindent
time computation (once TC and Galois groups computed): 16min, 47,728 ms\\
Curve: $17x^8 + (45y + 94)x^7 + (48y^2 + 218y + 249)x^6 + (26y^3 + 190y^2 + 520y + 422)x^5 + (6y^4 + 79y^3 + 384y^2 + 785y + 500)x^4 + (13y^4 + 131y^3 + 484y^2 + 785y + 422)x^3 + (18y^4 + 131y^3 + 384y^2 + 520y + 249)x^2 + (13y^4 + 79y^3 + 190y^2 + 218y + 94)x + (6y^4 + 26y^3 + 48y^2 + 45y + 17)=0$,                $P_0 = [0, -1, 1]$, $g = 5$\\
$\Sset = \{2, 67\}$, $\Tset = \{67, 9439, 277, 23, 293, 1163, 263, 17, 271, 239, 241\}$, $\Delta = \{ \textit{"small orbits"} \}$\\
Theta Fraction Rational TC: $67200yx^6 + (55440y^2 + 336000y)x^5 + (-5040y^3 + 154560y^2 + 739200y)x^4 + (-8400y^4 - 97440y^3 + 186480y^2 + 940800y)x^3 + (-42000y^4 - 184800y^3 + 186480y^2 + 739200y)x^2 + (-42000y^4 - 97440y^3 + 154560y^2 + 336000y)x + (-8400y^4 - 5040y^3 + 55440y^2 + 67200y)$\\
Orbit Decomposition. $\OTC: 1_1 15_1 30_2 60_1 180_2$, $\ETC: 3_1 45_1 60_2 90_2 180_1$\\
Orbit Decomposition. $J[2]: 3_1 15_1 30_2 45_1 60_3 90_2 180_3$\\
$\dim_{\F_2} J[2](\Q)= 0$\\
classgroup: 1,   1,   1,   1,   1\\   
$\dim_{\F_2} \kF(\Q)= 2$.  Global map injective? No\\
$\forall p \neq 2, \dim_{\F_2} \kF(\Q_p)=0$,      $\dim_{\F_2} \kF(\Q_2)=2$\\
$\dim_{\F_2} \Hop^1_{\Sset}(\Q,\E)= 74$\\
$\dim_{\F_2} \Sha^1(\Q,\E)=0$\\
$\dim_{\F_2} S_{\Tset}\leq 3$,	$\dim_{\F_2} \sKF\leq 2$,	$\dim_{\F_2} \Sel \leq 5$,	$\rank_\Q(J) \leq 5$\\
With careful look at $p=2$, $\rank_\Q(J) \leq 2$\\
rational points (not complete): $[0, -1, 1], [1, -3, 1], [-1, 1, 0]$\\
Bound by torsion-check: $1 \leq \rank_\Q(J)$ \hspace{5mm} $\leadsto 1 \leq \rank_\Q(J) \leq 2$\\
analytic rank: $2$

\newpage
\bibliography{bib}
\end{document}